\documentclass[a4paper,12pt]{amsproc}
\usepackage{amsmath,amsfonts,amssymb,amsthm,anysize}
\usepackage{enumerate}
\usepackage{epsfig}
\usepackage{supertabular}
\usepackage{graphicx}

\newcommand\Z{{\mathbb Z}}
\newcommand\Q{{\mathbb Q}}

\newcommand\C{{\mathbb C}}
\newcommand\F{{\mathbb F}}

\newcommand\Gal{{\mathrm{Gal}}}

\newcommand\Tr{{\mathrm{Tr}}}
\newcommand\Norm{{\mathrm{Norm}}}

\newcommand\ord{{\mathrm{ord}}}

\newcommand\lcm{{\mathrm{lcm}}}
\newcommand\inv{{\mathrm{inv}}}

\newcommand\Aut{\mathsf{Aut}}

\renewcommand\mod{{\mathrm{mod\, \, }}}

\newcommand\CC{\mathsf{D}}
\newcommand\D{\mathsf{D}}
\newcommand\QC{\mathsf{X}}

\newcommand\cP{{\mathcal P}}

\theoremstyle{plain}
\newtheorem{theorem}{Theorem}[section]

\newtheorem{lemma}[theorem]{Lemma}
\newtheorem{corollary}[theorem]{Corollary}

\newtheorem{proposition}[theorem]{Proposition}
\newtheorem{example}[theorem]{Example}

\numberwithin{equation}{section}
\allowdisplaybreaks[4]

\theoremstyle{remark}
\newtheorem{remark}[theorem]{Remark}

\usepackage{url, hyperref}
\hypersetup{citecolor=blue, linkcolor=blue, colorlinks=true}

\renewcommand\le{\leqslant}
\renewcommand\ge{\geqslant}

\begin{document}

\title[Pure Gauss sums]{Pure Gauss sums and skew Hadamard difference sets}

\author{Koji Momihara}
\address{ %
Division of Natural Science\\
Faculty of Advanced Science and Technology\\
Kumamoto University\\
2-40-1 Kurokami, Kumamoto 860-8555, Japan}
\email{momihara@educ.kumamoto-u.ac.jp}
\thanks{The author acknowledges the support by 
JSPS under Grant-in-Aid for Scientific Research (C) 20K03719.}

\subjclass[2010]{11L05, 11T22, 11T24, 05B10}
\keywords{}

\begin{abstract}
Chowla~(1962), McEliece~(1974), Evans~(1977, 1981) and Aoki~(1997, 2004, 2012)  studied Gauss sums, some integral powers of which are in the field of rational numbers. Such Gauss sums are called {\it pure}. 
In particular, Aoki (2004) gave a necessary and sufficient condition for a Gauss sum to be pure in terms of Dirichlet characters modulo the order of the multiplicative character involved. 
In this paper, we study pure Gauss sums with odd extension degree $f$ and classify them for $f=5,7,9,11,13,17,19,23$ based on Aoki's theorem.  Furthermore, we  characterize  a special subclass of pure Gauss sums in view of an application for skew Hadamard difference sets.  Based on the 
characterization, we give a new construction of skew Hadamard difference sets from cyclotomic classes of finite fields.  

%
\end{abstract}


\maketitle

\section{Introduction}\label{sec:intro}
Let $p$ be a prime and $f$ be a positive integer.  
Let $\F_{p^f}$ denote the finite field of order $p^f$. 
The canonical additive character $\psi$ of $\F_{p^f}$ is defined by 
$$\psi\colon\F_{p^f}\to \C^{\ast},\qquad\psi(x)=\zeta_p^{\Tr _{p^f/p}(x)},$$
where $\zeta_p={\rm exp}(\frac {2\pi i}{p})$ is a complex primitive $p$-th root of unity and $\Tr _{p^f/p}$ is the absolute trace from $\F_{p^f}$ to $\F_p$. 
All complex characters  of $(\F_{p^f},+)$ are given by $\psi_a$, where $a\in \F_{p^f}$. Here $\psi_a$ is defined by
\begin{equation}\label{additive}
\psi_a(x)=\psi(ax), \;\forall x\in \F_{p^f}. 
\end{equation}
Let $N$ be a positive divisor of $p^f-1$. 
For a multiplicative character 
$\eta_N$ of order $N$ of $\F_{p^f}$, we define the {\it Gauss sum} of $\F_{p^f}$ 
\[
G_{p^f}(\eta_N)=\sum_{x\in \F_{p^f}^\ast}\psi(x)\eta_N(x).
\] 

The Gauss sum is one of important and fundamental objects in number theory. 
The concept of Gauss sums was introduced by Gauss in 1801~\cite{Gauss}, who  evaluated the {\it quadratic} Gauss sums as in Theorem~\ref{thm:quad}. 
After Gauss' work, many researchers have tried to evaluate Gauss sums for larger $N$. 
However, in general, the explicit evaluation of Gauss sums is a very difficult problem. There are only a few cases where the Gauss sums have been completely evaluated. For example,  the Gauss sums for $N = 3, 4, 5, 6, 8, 12,16,24$ have been 
evaluated (but not explicit in some cases). %
See \cite{BEW97} for more details. 
The next important case is the so-called
{\it semi-primitive case} (also referred to as uniform cyclotomy or supersingular), where there
exists an integer $s$ such that $p^s \equiv -1 \,(\mod{N})$. See Theorem~\ref{thm:semiprim} for the explicit evaluation in this case.
The next interesting case is the {\it index $2$ case}, where the subgroup $\langle p\rangle$ generated by $p\in \Z$
has index $2$ in $(\Z/N\Z)^\times$. 
Many authors have studied  this case, see, e.g., \cite{Lang,M98,M03,YX09,YX10}.
In particular, a complete solution to the problem of evaluating Gauss sums in this case was given in \cite{YX10}. As a large generalization, Aoki~\cite{A10} studied  Gauss sums such that $(\Z/N\Z)^\times /\langle p\rangle$ is  an  elementary abelian $2$-group. The {\it index $4$ case} including the case where $(\Z/N\Z)^\times /\langle p\rangle$  is cyclic was also studied in \cite{FY,FYL,YLF}.

On the other hand, there were studies on Gauss sums from another point of view. 
Chowla~\cite{C01,C02} showed that if a Gauss sum defined in a prime field has the form $\epsilon p^{\frac{1}{2}}$ with $\epsilon$ a root of unity, it is in the quadratic case. 
McEliece~\cite{Mc} studied for which $(N,p,h)$, some nonzero integral power of the corresponding Gauss sum is an integer, i.e., $p^{-h/2}G_{p^h}(\eta_N)$ is a root of unity, related to weight distribution of irreducible cyclic codes.  Such Gauss sums are called {\it pure}. 
It is clear that the quadratic Gauss sums and the semi-primitive Gauss sums are examples of pure Gauss sums. 
Evans~\cite{E77} showed that  pure Gauss sums for prime powers $N$ are in the semi-primitive case. 
Furthermore,  Evans~\cite{E81} 
gave some nontrivial families of pure Gauss sums which are not semi-primitive. 
On the other hand, Aoki~\cite{A97} classified pure Gauss sums for small extension degrees  as follows.  
\begin{theorem}\label{thm:foddA}{\em \cite{A97}}
Assume that $f\in \{1,2,3,4\}$ and the order of $p$ modulo $N$ is $f$. Then, the corresponding Gauss sum $G_{p^f}(\eta_N)$ is pure if and only if it is of semi-primitive except for the following cases: 
\begin{align*}
f=3: &\,  (N,[p]_N)=(14,9),(14,11),(42,25),(42,37),(78,55),(78,71),\\ 
f=4: &\, (N,[p]_N)=(20,13),(20,17), (30,17),(30,23),(60,17),(60,53),(120,83),(120,107),
\end{align*}
where  $[p]_N$ is an integer such that $[p]_N\equiv p\,(\mod{N})$ and $1\le [p]_N\le N-1$. 
\end{theorem}
Furthermore, as a remarkable result, 
Aoki~\cite{A04,A12} gave a necessary and sufficient condition for a Gauss sum to be 
pure in terms of Dirichlet characters of modulo $N$, see Theorem~\ref{thm:aoki}. Based on the result,   
Aoki~\cite[Theorem~1.2]{A04} proved that 
for any  fixed $f$, the set of pairs $(N,[p]_N)$  such that 
the Gauss sum $G_{p^f}(\eta_N)$ is pure but not semi-primitive is finite. 

The evaluating Gauss sums is an important work also in view of applications in Combinatorics. 
In fact, Gauss sums have rich applications in the studies of combinatorial objects, such as difference sets, irreducible cyclic codes, strongly regular Cayley graphs,  cyclotomic association  schemes, sequences with good auto-correlation property,  highly nonlinear functions, etc. See, e.g., \cite{ADP,FWX,HH,Mc,MWX2019,S}. 
In particular,  pure Gauss sums were used 
for constructing skew Hadamard difference sets inequivalent to the classical Paley difference sets~\cite{CF15,FX113,Mo1}. 

Let $G$ be an additively written  group. We call a subset $D$ of $G$ a {\it difference set} if the list of differences 
``$x-y, x,y\in D,x\not=y$'' represents every element of $G\setminus \{0_G\}$ exactly $\lambda$ times. In this paper, we are concerned with difference sets in the additive group of the finite field, i.e., $G$ is an elementary abelian group. We say that a difference set is {\it skew Hadamard} if   $D$ is a skew-symmetric $(|G|-1)/2$-subset of $G$, i.e., $D\cup -D=G\setminus \{0_G\}$ and $D\cap -D=\emptyset$, where $-D=\{-x:x \in D\}$. 
The primary example of skew Hadamard difference sets is the classical {\it Paley difference set} in the additive group of the finite field $\F_q$ of order $q$ with $q\equiv 3\,(\mod{4})$, which consists of all nonzero squares of $\F_q$. 
The Paley difference set was the only known example in abelian groups for many years. Therefore, many researchers had believed that up to equivalence the Paley difference sets are the only skew Hadamard difference sets in elementary abelian groups. In 2006, Ding and Yuan~\cite{DY06} disproved this conjecture by  giving counterexamples of skew Hadamard difference sets in $(\F_{3^5},+)$. After their work, there have been many studies on constructions and classification of  skew Hadamard difference sets. See short surveys in Introduction of \cite{CF15,FX113,Mo1}. In particular,  Feng and Xiang~\cite{FX113} gave a construction of skew Hadamard difference sets based on pure Gauss sums, which are also in the index $2$ case. 
Furthermore, Chen and Feng~\cite{CF15} generalized the construction using pure Gauss sums satisfying $2\equiv p^j\,(\mod{N/2})$ for some integer $j$, 
see Theorem~\ref{thm:fx}.  
Their constructions are very flexible as explained in the next section, and give
rise to many skew Hadamard difference sets inequivalent to the Paley difference sets~\cite{Mo1}. 
The study in this paper is a continuation of those in \cite{CF15,FX113}.

In this paper, we will study pure Gauss sums with $f$ odd and their application for constructing skew Hadamard difference sets. 
The objectives of this paper are three-fold. First, we give some necessary conditions for pure Gauss sums with $f$ odd based on Aoki's Theorem~\ref{thm:aoki}, and update the result of Theorem~\ref{thm:foddA} for $f\in \{5,7,9,11,13,17,19,23\}$ in Theorem~\ref{thm:main111}. Second, we characterize pure Gauss sums such that $f$ is odd and $\langle p\rangle$ 
has index at most $8$ in $(\Z/N\Z)^\times$, and see that almost all pure Gauss sums  for $N\le 5000$ and odd $f$  fall into those classes. Third, we give a characterization for a special class of pure Gauss sums with the following property in view of applications for skew Hadamard difference sets: 
for $N=2m_1m_2\cdots m_r$,  
$G_{p^f}(\eta_{2\prod_{i\in J}m_i})$ is pure for any subset $J$ of $\{1,2,\ldots,r\}$ containing $1$, where $m_i$'s are distinct odd prime powers. 
Based on the characterization of pure Gauss sums, we give a new  construction of skew Hadamard difference sets from cyclotomic classes of  finite fields, which gives rise to two existence results. One of the results (that is, Corollary~\ref{cor:const1}) is covered by the result in \cite{CF15}, and the other (that is, Corollary~\ref{cor:const2}) is completely new not within the framework of previous studies. 
\section{Preliminaries}
\subsection{Basic properties of Gauss sums}
From the definition of Gauss sums, we see clearly that $G_{p^f}(\eta_N)$ is in the ring of algebraic integers of the field $\Q(\zeta_{p},\zeta_{N})$. Let $\sigma_{a,b}$ be the automorphism of $\Q(\zeta_{p},\zeta_{N})$ defined by 
\[
\sigma_{a,b}(\zeta_N)=\zeta_{N}^a, \qquad
\sigma_{a,b}(\zeta_p)=\zeta_{p}^b,
\]
where $\gcd{(a,N)}=\gcd{(b,p)}=1$. 
Below we list several basic properties of Gauss sums \cite{BEW97}. 
\begin{itemize}
\item[(i)] $G_{p^f}(\eta_N)\overline{G_{p^f}(\eta_N)}=p^f$ if $\eta_N$ is nontrivial. 
\item[(ii)] $G_{p^f}(\eta_N^p)=G_{p^f}(\eta_N)$.  
\item[(iii)] $G_{p^f}(\eta_N^{-1})=\eta_N(-1)\overline{G_{p^f}(\eta_N)}$.
\item[(iv)] $G_{p^f}(\eta_N)=-1$ if $\eta_N$ is trivial. 
\item[(v)] $\sigma_{a,b}(G_{p^f}(\eta_N))=\eta_N^{-a}(b)G_{p^f}(\eta_N^a)$.
\end{itemize}
In general, explicit evaluations of Gauss sums are very difficult. There are only a few cases where the Gauss sums have been evaluated. 
The most well-known case is the {\it quadratic} case, i.e., the $N=2$ case. 
\begin{theorem}\label{thm:quad}{\em (\cite[Theorem~5.15]{LN97})}
Let $\eta_2$ be the quadratic character of $\F_{p^f}$. Then, 
$G_{p^f}(\eta_2)=\epsilon p^{f/2}$, where 
\[
\epsilon=
\begin{cases}
(-1)^{f-1}&  \mbox{ if $p\equiv 1\,(\mod{4})$}, \\
(-1)^{f-1}i^f&  \mbox{ if $p\equiv 3\,(\mod{4})$}. 
\end{cases}
\]
\end{theorem}
The next simple case is the so-called {\it semi-primitive case} (also referred to as uniform cyclotomy or supersingular), where there
exists an integer $s$ such that $p^s \equiv -1 \,(\mod{N})$.  
\begin{theorem}\label{thm:semiprim}{\em (\cite[Theorem~11.6.1]{BEW97})} 
Suppose that $N>2$ and $p$ is semi-primitive modulo $N$, 
i.e., there exists an $s$ such that  $p^s\equiv -1\,(\mod{N})$. Choose 
$s$ minimal and write 
$h=2st$. Let $\eta_m$ be a multiplicative character of order $m$. 
Then, 
\[
p^{-h/2}G_{p^h}(\eta_m)=
\left\{
\begin{array}{ll}
(-1)^{t-1}&  \mbox{if $p=2$;}\\
(-1)^{t-1+(p^s+1)t/m}&  \mbox{if $p>2$. }
 \end{array}
\right.
\]
\end{theorem}

In this  paper, we will need the  {\it Davenport-Hasse product formula}, which is stated below.
\begin{theorem}
\label{thm:Stickel2}{\em (\cite[Theorem~11.3.5]{BEW97})}
Let $\theta$ be a multiplicative character of order $\ell>1$ of  $\F_{p^f}$. For  any nontrivial multiplicative character $\eta$ of $\F_{p^f}$, 
\begin{equation}\label{eq:HD}
G_{p^f}(\eta)=\frac{G_{p^f}(\eta^\ell)}{\eta^\ell(\ell)}
\prod_{i=1}^{\ell-1}
\frac{G_{p^f}(\theta^i)}{G_{p^f}(\eta\theta^i)}. 
\end{equation}
\end{theorem}
It is not easy to determine $\eta^\ell(\ell)$ in general. The following transformation of \eqref{eq:HD} is sometimes useful. 
\begin{corollary}\label{cor:HD}
With notation as in Theorem~\ref{thm:Stickel2}, if $\ell$ is odd, 
\[
G_{p^f}(\eta^\ell)=p^{-f\frac{\ell-1}{2}}\sigma_{1,\ell^{-1}}
\Big(\prod_{i=0}^{\ell-1}G_{p^f}(\eta\theta^i)\Big).
\] 
\end{corollary}
\proof
Note that
\[
\prod_{i=1}^{\ell-1}G_{p^f}(\theta^i)=\prod_{i=1}^{\frac{\ell-1}{2}}G_{p^f}(\theta^i)G_{p^f}(\theta^{\ell-i})=p^{f\frac{\ell-1}{2}}\prod_{i=1}^{\frac{\ell-1}{2}}\theta^i(-1). 
\]
Here, $\theta(-1)=1$; otherwise $\theta^{\ell}(-1)=-1$, a contradiction to that $\theta^{\ell}$ is trivial.   
Hence, 
$\prod_{i=1}^{\ell-1}G_{p^f}(\theta^i)=p^{f\frac{\ell-1}{2}}$. Furthermore, 
note that \[
G_{p^f}(\eta^\ell)\eta^{-\ell}(\ell)=\sigma_{1,\ell}(G_{p^f}(\eta^\ell)). 
\]
Then, \eqref{eq:HD} is reformulated as 
\begin{equation}\label{eq:HDp}
\sigma_{1,\ell}(G_{p^f}(\eta^\ell))=p^{-f\frac{\ell-1}{2}}\prod_{i=0}^{\ell-1}G_{p^f}(\eta\theta^i). 
\end{equation}
Finally, by acting $\sigma_{1,\ell^{-1}}$ to both sides of \eqref{eq:HDp}, we obtain the assertion of the corollary. 
\qed
\vspace{0.3cm}

We will also need the  {\it Davenport-Haase lifting formula}, which is stated below.
\begin{theorem}\label{thm:lift}{\em (\cite[Theorem~5.14]{LN97})}
Let $\eta$ be a nontrivial multiplicative character of $\F_{p^f}$ and 
let $\eta'$ be the lift of $\eta$ to $\F_{p^{fs}}$, i.e., $\eta'(\alpha)=\eta(\Norm_{p^{fs}/p^f}(\alpha))$ for $\alpha\in \F_{p^{fs}}$, where $s\geq 2$ is an integer. Then 
\[
G_{p^{fs}}(\eta')=(-1)^{s-1}(G_{p^f}(\eta))^s. 
\]
\end{theorem}

\subsection{Pure Gauss sums}\label{sec:puredef}
Let $\eta_N$ be a multiplicative character of order $N$ of $\F_{p^f}$. We say that the Gauss sum $G_{p^f}(\eta_N)$ is {\it pure} if $\epsilon=G_{p^f}(\eta_N)p^{-f/2}$ is a root of unity. We call the $\epsilon$ as the {\it sign} or the {\it root of unity} of the pure Gauss sum $G_{p^f}(\eta_N)$.  
\begin{lemma}
If $G_{p^f}(\eta_N)$ is pure, so is $G_{p^f}(\eta_N^a)$ for any  $a$ with  $\gcd{(a,N)}=1$. 
\end{lemma}
\proof
Let $\sigma_{a,1}\in \Gal(\Q(\zeta_{p},\zeta_{N})/\Q)$. Then, $\sigma_{a,1}(G_{p^f}(\eta_N))=G_{p^f}(\eta_N^a)$ is also pure. 
\qed

\vspace{0.3cm}
The lemma above implies that the purity of Gauss sums is depending on $N$ but not depending on the choice of $\eta_N$. Then, denote by  ${\mathcal P}$ the set of triples $(N,f,p)$ such that $G_{p^f}(\eta_N)$ is pure. 

From now on, let $f$ be the order of $p$ modulo $N$ and  $\eta_N$ be a multiplicative character of order $N$ of $\F_{p^f}$. Let $s$ be any positive integer and $\eta_N'$ be the lift of $\eta_N$ to $\F_{p^{fs}}$. If $G_{p^f}(\eta_N)$ is pure, then 
so is $G_{p^{fs}}(\eta_N')$ by Theorem~\ref{thm:lift}. Hence, the purity problem of $G_{p^{fs}}(\eta_N')$ is reduced to that of $G_{p^f}(\eta_N)$. 
Hence, we 
consider 
\[
{\mathcal P}^\ast:=\{(N,f,p)\in {\mathcal P}\mid \ord_{N}(p)=f\}. 
\]
 
The following characterization of pure Gauss sums is obtained from the well-known {\it Stickelberger} theorem on ideal factorizations of Gauss sums~\cite[Theorem~11.2.2]{BEW97}. 
\begin{proposition}\label{prop:aoki0}{\em (\cite{A04,E77,KL79})}
$(N,f,p)\in {\mathcal P}^\ast$ if and only if 
\[
\sum_{i=0}^{f-1}[tp^i]_N=\frac{fN}{2}
\]
for any integer $t$ prime to $N$, where $[x]_N$ is an integer such that $0\le [x]_N\le N-1$ and $[x]_N\equiv x\,(\mod{N})$.  
\end{proposition}
The proposition above gives the following characterization.  
\begin{lemma}{\em (\cite{E77,E81})}\label{lem:evcha}
If $(N,f,p)\in {\mathcal P}^\ast$, it holds that $N\,|\,(p^f-1)/(p-1)$ or  
$N/2\,|\,(p^f-1)/(p-1)$ depending on whether $f$ is even or odd. 
\end{lemma}
On the other hand, Proposition~\ref{prop:aoki0} implies that the purity of  Gauss sums for a fixed $N$ depends only on the residue class of $p$ modulo $N$. Furthermore, 
the proposition implies that if $(N,f,p) \in {\mathcal P}^\ast$, then $(N,f,r)\in {\mathcal P}^\ast$ for any prime $r\equiv p^i\,(\mod{N})$, where $i$ is an arbitrary integer such that $1\le i\le f-1$ and $\gcd{(i,f)}=1$.   

The Gauss sums in semi-primitive case are clearly pure. 
Hence, we have  
\begin{align*}
({\mathcal P}^{(-1)}:=)\{(N,f,p)\mid \exists i \mbox{ s.t. }p^i\equiv -1\,(\mod{N})\}\subseteq 
{\mathcal P}. 
\end{align*}
It is clear that 
$f$ is even if $(N,f,p)\in {\mathcal P}^{(-1)}$. 
Evans~\cite{E77} showed that  pure Gauss sums for prime powers $m$ are in the semi-primitive case. 
Furthermore, Evans~\cite{E81} also gave the following nontrivial sufficient conditions for Gauss sums to be pure. 
\begin{theorem}\label{thm:small}
Suppose that $m=cd$ with $\gcd{(c,d)}=\gcd{(\ord_{c}(p),\ord_{d}(p))}=1$ and let $f=\ord_{m}(p)$, where $\ord_{n}(x)$ is the order of $x$ in $(\Z/n\Z)^\times$. Then, $G_{p^f}(\eta_m)$ is pure if any of the following holds. 
\begin{itemize}
\item[(1)] $\ord_{c}(p)=\phi(c)$ and $\ell\in \langle p\rangle\,(\mod{d})$ for some prime $\ell\,|\,c$. 
\item[(2)] $-1\not\in \langle p\rangle\,(\mod{c})$, $2\ord_{c}(p)=\phi(c)$, $\ell\in \langle p\rangle\,(\mod{d})$ for some prime $\ell\,|\,c$, and all of them hold with $c$ and $d$ interchanged. 
\item[(3)] $2||m$, $2+m/2\not\in \langle p\rangle\,(\mod{c})$,  $2\ord_{c}(p)=\phi(c)$, $-1$ or $\ell$ is in $\langle p\rangle\,(\mod{d})$ for some prime $\ell\,|\,c$, and all of them hold with $c$ and $d$ interchanged. 
\end{itemize} 
Here, $\phi$ is Euler's totient function.  
\end{theorem}
On the other hand, Aoki~\cite[Theorem~7.2]{A04} proved that the converse of the assertion of Theorem~\ref{thm:small} also  holds if $c$ and $d$ are both odd prime powers. 

In this paper, we are concerned with pure Gauss sums with $f$ odd. 
There were not so many studies on pure Gauss sums for odd $f$ in the literature. 
\begin{proposition}\label{prop:111}{\em (\cite[Corollary~3]{E81})}
If $f$ is odd and $(N,f,p)\in {\mathcal P}^\ast$, then $2\|N$. 
\end{proposition}

The following proposition comes from Theorem~\ref{thm:small}~(1) as   
$c=\ell=2$ or Corollary~8 in \cite{E81}. 
(Note that Theorem~\ref{thm:aoki} below is a large generalization of Theorem~\ref{thm:small}.) 
Chen-Feng~\cite{CF15} also gave a proof for the result based on Davenport-Hasse product formula. 
\begin{proposition}\label{prop:2m}
Assume that $2\|N$. If there exists $j$ such that $p^j\equiv 2\,(\mod{N/2})$, then $(N,f,p)\in {\mathcal P}^\ast$. In particular, $G_{p^f}(\eta_N)=G_{p^f}(\eta_2)$. 
\end{proposition}
The proposition above defines a class of pure Gauss sums with $f$ odd: 
\[
{\mathcal P}^{(2)}:=\{(N,f,p)\mid f \mbox{ is odd, }\exists i \mbox{ s.t. }p^{i}\equiv 2\,(\mod{N/2})\}\subseteq 
{\mathcal P}. 
\] 
We will study in Section~\ref{sec:appli} whether there is a class of pure Gauss sums other than ${\mathcal P}^{(2)}$ compatible with a construction of skew Hadamard difference sets. 

Next, we give one basic property of pure Gauss sums. 
\begin{lemma}\label{prop:2mR}
Let $p_1$ be an odd prime and $t$ be a positive integer.  
Assume that $p_1^t\|N$.  
If $(N,f,p)\in {\mathcal P}$, then  $(N/p_1^s,f,p)\in {\mathcal P}$ for any $s\le t-1$. 
\end{lemma}
\proof
Let $\eta_N$ be a multiplicative character of order $N$ of $\F_{p^f}$ and 
$\theta$ be a multiplicative character of order $p_1^s$ of $\F_{p^f}$. 
Then, by Theorem~\ref{thm:Stickel2} as $\ell=p_1^s$, we have  
\[
G_{p^f}(\eta_N)=\frac{G_{p^f}(\eta_N^\ell)}{\eta_N^\ell(\ell)}
\prod_{i=1}^{\ell-1}
\frac{G_{p^f}(\theta^i)}{G_{p^f}(\eta_N\theta^i)}.  
\]
Since $\eta_N\theta^i$ is of order $N$, $G_{p^f}(\eta_N\theta^i)$ is also pure. 
On the other hand, $\prod_{i=1}^{\ell-1}
G_{p^f}(\theta^i)=\Big(\prod_{i=1}^{\frac{\ell-1}{2}}\theta^i(-1)\Big) p^{f\frac{\ell-1}{2}}$. Hence, 
$G_{p^f}(\eta_N^\ell)$ is pure, i.e., $(N/p_1^s,f,p)\in {\mathcal P}$. 
 \qed

\vspace{0.3cm}
We will need the following powerful characterization of pure Gauss sums given by Aoki~\cite{A12}. To state it, 
let $\D(N)$ denote the set of Dirichlet characters modulo $N$, and 
define 
\begin{align*}
\CC^-(N,p):=&\,\{\chi\in \D(N)\mid \chi(p)=1,\chi \mbox{ is an odd character}\}, \\
\QC^-(N,p):=&\,\{\chi\in \CC^-(N,p)\mid \mbox{The conductor of $\chi$ is divisible by any prime factor of $N$}\}. 
\end{align*}
\begin{theorem}\label{thm:aoki}
{\em (\cite[Proposition 4.9]{A12})}
$(N,f,p)\in {\mathcal P}^\ast$ if and only if 
the following two conditions hold. 
\begin{itemize}
\item[(1)] $\QC^{-}(N,p)=\emptyset$. 
\item[(2)] For 
any $\chi\in \D^-(N,p)$, there exists a prime divisor $\ell$ of $N$ but not dividing the conductor of $\chi$ such that $\chi(\ell)=1$. 
\end{itemize}
\end{theorem}
Define 
\[
{\mathcal P}^\ast_f:=\{(N,\overline{p})\mid (N,f,p)\in {\mathcal P}^\ast\setminus {\mathcal P}^{(-1)}\}, 
\]
where $\overline{p}$ denotes a minimum representative in $\{[p^i]_N\mid 1\le i\le f-1,\gcd{(i,f)}=1\}$. 
Based on the theorem above, Aoki~\cite[Theorem~1.1]{A12} proved that ${\mathcal P}^\ast_f$ is a finite set for every positive integer $f$. 

\subsection{Skew Hadamard difference sets}
Let $q=p^f$ be a prime power  and let $N>1$ be a divisor of $q-1$. Let $C_i^{(N,q)}=\omega^i \langle \omega^N\rangle$, $0\le i\le N-1$, be 
the {\it cyclotomic classes} of order $N$ of $\F_q$, where $\omega$ is a fixed primitive element of $\F_q$. We assume that $q\equiv 3\,(\mod{4})$ and $N$ is even. Then, it is clear that $N/2$ is odd. 
In this paper, we will give a construction for a skew Hadamard difference set $D$ as a union of suitable  $m=N/2$ cyclotomic classes. To do this, we will use the following well-known characterization of skew Hadamard difference sets.  
\begin{lemma}\label{lem:bas}
Let $G$ be an abelian group of order $v\equiv 3\,(\mod{4})$, and let $D$ be a skew symmetric $(v-1)/2$-subset of $G$.  
The set $D$ is a skew Hadamard difference set if and only if $\psi(D)\in \{\frac{-1+\sqrt{-v}}{2},\frac{-1-\sqrt{-v}}{2}\}$ for any 
nontrivial character of $G$. 
\end{lemma}
Note that $D^\perp=\{\psi\in G^\perp\,|\,\psi(D)=\frac{-1+ \sqrt{-q}}{2}\}$  (and its inverse) also forms a skew Hadamard difference set, called the {\it dual} of  $D$, in the character group $G^\perp$ of $G$. 
The following result is also known (cf. \cite[Lemma~2.1]{CXS}). 
\begin{lemma}\label{lem:skewcheck}
Let $G$ be an abelian group of order $p^h$, where $p$ is a prime such that $p\equiv 3\,(\mod{4})$ and $h$ is an odd integer. Let $D$ be a skew symmetric $(p^h-1)/2$-subset of $G$ such that $D$ is invariant under the multiplication by $x^2$ for $x\in \F_p^\ast$. If $2\psi(D)+1\equiv 0\,(\mod{p^{\frac{h-1}{2}}})$ for any nontrivial character $\psi$ of $G$,  then $D$ is a skew Hadamard difference set in $G$.  
\end{lemma}

Let $I$ be a $N/2$-subset of $\{0,1,\ldots,N-1\}$. 
To check whether a candidate subset $D=\bigcup_{i\in I}C_i^{(N,q)}$ is a skew Hadamard difference set, by Lemma~\ref{lem:bas}, it suffices to show that $(\psi_a(D):=)\sum_{x\in D}\psi_a(x)\in \{\frac{-1\pm \sqrt{-q}}{2}\}$ for any nonzero $a\in \F_q$. 
Note that the character value $\psi_a(D)$ can be expressed as a linear combination of 
Gauss sums  (cf.~\cite{FX113}) by using the orthogonality of characters: 
\begin{eqnarray}\label{eigen}
\psi_a(D)=\frac{1}{N}
\sum_{i=0}^{N-1}G_{q}(\eta_N^{-i})
\sum_{i\in I}\eta_N(a\gamma^i ), 
\end{eqnarray}
where $\eta_N$ is a fixed multiplicative character of order $N$ of $\F_q$. 
Thus, the computations needed to check whether a candidate subset $D=\bigcup_{i\in I}C_i^{(N,q)}$ is a skew Hadamard difference set are essentially reduced to evaluating Gauss sums. 
For example, if $N=2$, 
we have 
\begin{equation}\label{eq:Gaussquad}
\psi_a(C_i^{(2,q)})=\frac{-1+(-1)^{a+i}G_q(\eta_2)}{2},  
\end{equation}
where $\eta_2$ is the quadratic character of $\F_q$. By Theorem~\ref{thm:quad}, we have  $\psi_a(C_i^{(2,q)})\in \{\frac{-1\pm \sqrt{-q}}{2}\}$ if $q\equiv 3\,(\mod{4})$. Hence, each $C_i^{(2,q)}$, $i=0,1$, is a skew Hadamard difference set in $(\F_q,+)$, that is, the so-called Paley difference set.


Let $X$ be a subset of $\F_{q^\ell}^\ast/\F_{q}^\ast$, and $\pi:\F_{q^\ell}^\ast\to \F_{q^\ell}^\ast/\F_{q}^\ast$ be the natural projection 
homomorphism. Define 
\[
D(X)=\{x\in C_0^{(2,q^\ell)}\,|\,\pi(x)\in X\}\cup \{x\in C_1^{(2,q^\ell)}\,|\,\pi(x)\not\in X\}.  
\] 
Chen-Feng~\cite{CF15} showed that under the assumptions that $\ell$ is odd and $X$ is a difference set with parameters $((q^\ell-1)/(q-1),q^{\ell-1},q^{\ell-2}(q-1))$, $D(X)$ is a skew Hadamard difference set or a Paley type partial difference set if and only if  $X$ is an {\it Arasu-Dillon-Player difference set}. Furthermore, they gave the following construction of skew Hadamard difference sets based on the class ${\mathcal P}^{(2)}$ of pure Gauss sums, which is a generalization of that given by Feng-Xiang~\cite{FX113}. 
\begin{theorem}\label{thm:fx}
{\em (\cite[Theorem~1.4]{CF15})} Let $q=p^f\equiv 3\,(\mod{4})$ be a prime power with $p$ a prime, and let $\ell$ be any odd positive integer. Let $m$ be a divisor of $(q^\ell-1)/(q-1)$ satisfying $2\equiv p^j\,(\mod{m})$ for some integer $j$, and $\tau:\F_{q^\ell}^\ast/\F_{q}^\ast\to \Z/m\Z$ be the natural projection. Then,  for any subset $X$ of $\Z/m\Z$, the set 
$D(\tau^{-1}(X))$ is  a skew Hadamard difference set in $(\F_{q^{\ell}},+)$. 
\end{theorem}
The original statement of the theorem above in \cite{CF15} assumed that $q^\ell$ is an odd prime power not necessarily $q^\ell\equiv 3\,(\mod{4})$ and  $m$ is a divisor of  $(q^\ell-1)/(q-1)$ satisfying  $-1\equiv p^j\,(\mod{m})$ or $2\equiv p^j\,(\mod{m})$  for some integer $j$ 
since the authors treated also Paley type partial difference sets not only skew Hadamard difference sets. However,  in our situation,  $-1\equiv p^j\,(\mod{m})$ is impossible since $f\ell$ is odd.  
Note that the assumption $2\equiv p^j\,(\mod{m})$ implies that $(p^j-1)\equiv 1\,(\mod{m})$ for some positive integer $j$. Hence, $p-1$ and $m$ are coprime, i.e., 
$m$ is a divisor of $(p^f-1)/(p-1)$. Hence, we can not remove the condition $m\,|\,(q^\ell-1)/(q-1)$. 
It is clear that $D(\tau^{-1}(X))$ in the theorem above is a union of cyclotomic classes of order $N=2m$ of $\F_{q^\ell}$. In particular, it is  expressed as 
$D(\tau^{-1}(X))=\bigcup_{i\in I}C_i^{(2m,q^\ell)}$, where $I=\{(m+1)i\,(\mod{2m})\,|\,i\in X\} \cup \{(m+1)i+m\,(\mod{2m})\,|\,i\in (\Z/m\Z)\setminus X\}$. In other words, one can take $I$ as an arbitrary $m$-subset of $\Z/2m\Z$ such that 
$\{i\,(\mod{m})\,|\,i\in I\}=\{0,1,\ldots,m-1\}$. 
Thus, the theorem above is very powerful. In fact, the choice of the set $X$ is very flexible and $f\ell$ can be taken as an arbitrary odd positive integer divisible by  the order of $p$ in $(\Z/m\Z)^\times$. 
Then,  Theorem~\ref{thm:fx} yields infinite families of skew Hadamard difference sets inequivalent to the Paley difference sets~\cite{Mo1}. 
In Section~\ref{sec:appli}, we give a construction of skew Hadamard difference sets based on a class of pure Gauss sums 
not belonging to ${\mathcal P}^{(2)}$. 

\section{Some necessary conditions for pure  Gauss sums with $f$ odd}

Throughout this section, we assume that $f=\ord_{N}(p)$ is odd  and $(N,f,p)\in {\mathcal P}^\ast$. Then, by Proposition~\ref{prop:111}, we have $2\| N$. 
Let $p_i$, $i=1,2,\ldots,r$, be distinct odd primes and $u_i$, $i=1,2,\ldots,r$, be
positive integers. 
Let $N=2m=2m_1m_2\cdots m_r$, where $m_i=p_i^{u_i}$, $i=1,2,\ldots,r$. 
Let $f_i$, $i=1,2,\ldots,r$, denote the orders of $p$ modulo $m_i$, respectively. Then, $f=\lcm{(f_1,f_2,\ldots,f_r)}$. 

In this section, we characterize $(N,\overline{p})\in {\mathcal P}^\ast_f$  for $f\in \{3,5,7,9,11,13,17,19,23\}$,  $N\le 5000$ or $\phi(N)/f\le 8$ with $f$ odd.   
\subsection{Necessary conditions}\label{subsec:nece}
Aoki~\cite{A04} proved the following theorem. 
\begin{theorem}{\it \cite[Theorem~5.1]{A04}}\label{thm:ao1}
Assume that $f$ is odd and $(N,f,p)\in {\mathcal P}^\ast$. 
\begin{itemize}
\item[(1)] If $r$ is odd,  $m_i\,|\,2^{2f}-1$ for each $i$, $i=1,2,\ldots,r$; 
\item[(2)]  If $r$ is even,  either $m_i\,|\,2^{2f}-1$ or $\phi(m_i)\,|\,4f$  for each $i$, $i=1,2,\ldots,r$.  
\end{itemize}
\end{theorem}
Since $\phi(m_i)\,|\,4f$ implies that  $m_i\,|\,2^{4f}-1$ by Fermat's little theorem, we have the following corollary. 
\begin{corollary}\label{cor:ao1}{\em (\cite[Corollary~5.2]{A04})}
Assume that $f$ is odd and $(N,f,p)\in {\mathcal P}^\ast$.  Then, $m\,|\,2^{4f}-1$. 
\end{corollary}
The corollary above implies that the set of pairs $(N,\overline{p})\in {\mathcal P}^\ast_f$ is finite for any fixed $f$. 
To classify $(N=2m,\overline{p})\in {\mathcal P}^\ast_f$ for a fixed odd $f$, we may take positive divisors $m$ of  $2^{4f}-1$ in view of Corollary~\ref{cor:ao1}. However, even if  $f$ is small,  some divisor $m$ of $2^{4f}-1$ is too large to check whether $(N,\overline{p})\in {\mathcal P}^\ast_f$ by a computer.   So, we will give some new necessary conditions for divisors $m$ of $2^{4f}-1$ such that $(N,\overline{p})\in {\mathcal P}^\ast_f$, which are all based on Aoki's Theorem~\ref{thm:aoki}. 

Let $\chi_i$ be a character of order $\phi(m_i)$ of $(\Z/m_i\Z)^\times$. Then, we have 
$\chi_i^{f_i}(p)=\chi_i(p^{f_i})=1$. Furthermore,  $\chi_i^{f_i}(-1)=\chi_i(-1)=-1$ since $f_i$ is odd. Hence, $\chi_i^{f_i}\in \D^-(N,p)$ for  $i=1,2,\ldots,r$. 
\begin{proposition}\label{nece0}
Assume that $f$ is odd and  $(N,f,p)\in {\mathcal P}^\ast$. 
If there is $j$ such that $m_j \not | 2^{f}-1$, there is $h$ with $h\not=j$ such that  $m_j\,|\,p_h^{f}-1$. 
\end{proposition}
\proof
By Theorem~\ref{thm:aoki}, we have 
$\chi_j^{f_j}(2)=1$ or $\chi_j^{f_j}(p_h)=1$ for some $h=1,2,\ldots,r$ with $h\not=j$. If $\chi_j^{f_j}(2)=\chi_j(2^{f_j})=1$, we have $m_j\,|\,2^{f_j}-1$, which contradicts to $m_j\not|2^{f}-1$. Hence, we have 
 $\chi_j^{f_j}(p_h)=1$ for some $h=1,2,\ldots,r$ with $h\not=j$. This implies that $m_j\,|\,p_h^{f_j}-1\,|\,p_h^{f}-1$. 
\qed
\vspace{0.3cm}

Next, we give two necessary conditions  for $(N,f,p)\in {\mathcal P}^\ast$ with $r$ even. 
\begin{proposition}\label{nece1}
Assume that $f$ is odd, $r$ is even and  $(N,f,p)\in {\mathcal P}^\ast$. 
If there is $j$ such that $m_j\not |2^{2f}-1$, it holds that $\phi(m_h)\le 2f$ for any $m_h$ such that $m_h\,|\,2^{2f}-1$. 
\end{proposition}
\proof 
Assume that $\phi(m_h)>2f$ for some $h$ such that $m_h\,|\,2^{2f}-1$. 
Note that $\chi_h^{2f}$ is nontrivial. Consider the character 
\[
\theta=\prod_{i=1}^r\chi_i^{f_i}. 
\]  
Since $r$ is even, $\theta$ is an even character. 
Next, we consider the characters 
\[
\theta'=\chi_{h}^{2f-f_h}\theta=\chi_{h}^{2f}\prod_{i\not=h}\chi_i^{f_i}
\]
and 
\[
\theta''=\prod_{i\not=h}\chi_i^{f_i}. 
\]
Since $\theta'\in \QC^-(m,p)$, by Theorem~\ref{thm:aoki}, we have  $\theta'(2)=1$. On the other hand, since $\chi_{h}^{2f}(2)=\chi_{h}(2^{2f})=1$,  
we have $\theta''(2)=1$. Similarly, for any $k$ with $k\not=h$, let 
\[
\theta'''=\chi_k^{-f_k}\prod_{i\not=h,k}\chi_i^{f_i}. 
\]
Then, $\theta'''(2)=1$.
Hence, we have $\chi_k^{2f_k}(2)=1$. Then,  
$m_k\,|\,2^{2f}-1$ for any $k$ with $k\not=h$, which contradicts to  $m_j\not|2^{2f}-1$ for some $j$. Hence, $\phi(m_h)\le  2f$ for any $h$ such that 
$m_h\,|\,2^{2f}-1$. 
\qed

\begin{remark}\label{rem:nece1}
We can improve Proposition~\ref{nece1} in the $r=2$ case as 
``If $m_1\not |2^{f}-1$ and $m_2\,|\,2^{2f}-1$, it holds that 
$\phi(m_2)\le 2f$.''  
Let $j=1$ and $h=2$ in the proof of Proposition~\ref{nece1}. 
Then, we can similarly prove that $\theta''(2)=\chi_1^{f_1}(2)=1$. Then, we have $m_1\,|\,2^{f_1}-1$, which contradicts to $m_1\not|2^f-1$. 
\end{remark}
\begin{proposition}\label{nece2}
Assume that $f$ is odd, $r$ is even and  $(N,f,p)\in {\mathcal P}^\ast$. 
If there is $j$ such that $m_j \not|2^{2f}-1$,  
it holds that $m_k\,|\,p_j^{2f}-1$ for any $k$ with $k\not=j$. 
\end{proposition}
\proof 
If $\phi(m_j)\le 2f_j$, since 
$2f_j$ divides $\phi(m_j)$, we have $\phi(m_j)=2f_j$, which implies that 
$\phi(m_j)\,|\,2f$. However, this contradicts to $m_j \not|2^{2f}-1$. Hence, we have $\phi(m_j)> 2f_j$, and then 
 $\chi_j^{2f_j}$ is nontrivial. 
Let 
\[
\theta=\chi_j^{f_j}\prod_{i=1}^r\chi_i^{f_i}=\chi_j^{2f_j}\prod_{i\not=j}\chi_i^{f_i}. 
\]  
Since  $\theta \in \QC^-(m,p)$, we have $\theta(2)=1$. Let 
\[
\theta'=\prod_{i\not=j}\chi_i^{f_i}.  
\]
Since $\theta'\in \QC^-(m/m_j,p)$, by Theorem~\ref{thm:aoki}, we have  either $\theta'(2)=1$ or $\theta'(p_j)=1$. 
If $\theta'(2)=1$, we have $\chi_j^{2f_j}(2)=1$, which contradicts to 
$m_j\not |2^{2f}-1$. Hence, 
$\theta'(p_j)=1$. For any $m_k$ with $k\not=j$, let  
\[
\theta''=\chi_{k}^{-f_k}\prod_{i\not=j,k}\chi_i^{f_i}.  
\]
Then, we similarly have $\theta''(p_j)=1$. Hence, $\chi_{k}^{2f_k}(p_j)=1$. This implies that 
$m_k\,|\,p_j^{2f}-1$. 
\qed
\vspace{0.3cm}


The statement of the proposition above is similar to \cite[Theorem~11.1]{A04} but not exactly same. 
We next give a necessary condition  for $(N,f,p)\in {\mathcal P}^\ast$ with $r$ odd. 

\begin{proposition}\label{nece3}
Assume that $f$ is odd, $r$ is odd and  $(N,f,p)\in {\mathcal P}^\ast$. 
If there are $j$ and $h$ such that $m_j\not |2^{f}-1$ and $m_h\,|\,2^{f}-1$,   then 
 either $\phi(m_h)\le 2f$ or $m_k\,|\,p_j^{2f}-1$ for any $k$ with $k\not=j,h$ and $m_h\,|\,p_j^{4f}- 1$. 
\end{proposition}
\proof 
If $\chi_j^{f_j}(2)=1$, we have $m_j\,|\,2^{f_j}-1$, which contradicts to $m_j\not|2^{f}-1$. Hence,  $\chi_j^{f_j}(2)\not=1$. 

Assume that $\phi(m_h)> 2f$. Noting that $\chi_h^{2f}$ is nontrivial, let 
\[
\theta=\chi_h^f \prod_{i\not=h}\chi_i^{f_i}
\]  
and 
\[
\theta'=\chi_h^{2f}\prod_{i\not=j,h}\chi_i^{f_i}. 
\]
Since $\theta  \in \QC^-(m,p)$, we have $\theta(2)=1$ by Theorem~\ref{thm:aoki}. 
Furthermore, since $\theta' \in \QC^-(m/m_j,p)$, we have  $\theta'(2)=1$ or $\theta'(p_j)=1$. 
If  $\theta'(2)=1$, we have $\chi_h^{f}(2)=\chi_j^{f_j}(2)\not=1$, which  
contradicts to $m_h\,|\,2^{f}-1$. Hence, $\theta'(p_j)=1$. For $k$ with $k\not=j,h$, 
let \[
\theta''=\chi_h^{2f}\chi_k^{-f_k}\prod_{i\not=j,h,k}\chi_i^{f_i}. 
\]
Then, we similarly have $\theta''(p_j)=1$. Hence, we obtain 
$\chi_k^{2f_k}(p_j)=1$. This implies that $m_k\,|\,p_j^{2f}- 1$. 
Furthermore, let 
\[
\theta'''=\chi_h^{-2f}\prod_{i\not=j,h}\chi_i^{f_i}. 
\]
Then, we similarly have $\theta'''(p_j)=1$. Hence, we obtain 
$\chi_h^{4f}(p_j)=1$. This implies that $m_h\,|\,p_j^{4f}-1$. 
 \qed

\begin{example}\label{exam:7}
Let $f=7$. Then, all integers $N>2$ satisfying the condition of Theorem~\ref{thm:ao1} are 
\begin{align*}
& 6, 30, 86, 174,254, 258, 290, 430, 762, 1270, 2494, 7366, \\
&10922,32766,   37410, 110490, 163830, 950214, 1583690.
\end{align*}
Proposition~\ref{nece0} reduces the list above to 
$ 254 , 762, 10922 , 32766$. Furthermore, $762$ and $10922 $ are 
excluded by Remark~\ref{rem:nece1}, and $32766$ is excluded by Proposition~\ref{nece3}. For the  remaining $N=254$, we have $(254,129)\in {\mathcal P}_7^\ast$. Thus,  $(N,\overline{p})\in {\mathcal P}^\ast_7$ are classified.   
\end{example}
The following is our main theorem in this subsection. 
\begin{theorem}\label{thm:main111}
For $f=3,5,7,9,11,13,17,19,23$, all $(N,\overline{p})\in {\mathcal P}^\ast_f$ are listed below: 
\begin{align*}
f=3: &\quad (N,\overline{p})=(14,9),(42,25),(78,55);\\
f=5: &\quad (N,\overline{p})=(62,33),(110,31);\\
f=7: &\quad (N,\overline{p})=(254,129);\\
f=9: &\quad (N,\overline{p})=(146,37),(1022,513); \\
f=11: &\quad (N,\overline{p})=(46,3),(178,39),(4094,2049); \\
f=13: &\quad (N,\overline{p})=(16382,8193); \\
f=17: &\quad (N,\overline{p})=(262142,131073);\\ 
f=19: &\quad (N,\overline{p})=(1048574,524289); \\
f=23: &\quad (N,\overline{p})=(94,3),(356962,83663), (16777214,8388609). 
\end{align*} 
\end{theorem}
\proof 
First, we list all $N=2m$ satisfying the condition of Theorem~\ref{thm:ao1}. 
Then, similarly to Example~\ref{exam:7},  we reduce the candidates of $N$ such that $(N,f,p)\in {\mathcal P}^\ast$ by applying Propositions~\ref{nece0}, \ref{nece1}, \ref{nece2}, \ref{nece3}, and  Remark~\ref{rem:nece1}.  For remaining candidates, we  used a computer to  directly check whether there is $p$ such that $(N,f,p)\in {\mathcal P}^\ast$ based on Proposition~\ref{prop:aoki0}.  
\qed
\subsection{Characterization of pure Gauss sums of small index}\label{subsec:chara}
In Tables~\ref{Tab0} and \ref{Tab01} of  the appendix, we will give a list of   
$(N,\overline{p})\in {\mathcal P}^\ast_f$ for $N\le 5000$ and odd $f$ by using a computer. Almost all examples listed in the tables belong to ${\mathcal P}^{(2)}$ or satisfy $\phi(N)/f\le 8$. 
Therefore, in this subsection, we characterize $(N,\overline{p})\in {\mathcal P}^\ast_f$ such that $f$ is odd and  $\phi(N)/f\le 8$.  
Note that $\phi(N)/f$ must be even since $f$ is odd. Hence, we consider the cases where $\phi(N)/f=2,4,6,8$. 
\begin{proposition} \label{prop:22}
Assume that $\phi(N)/f=2$. 
Then, $(N,f,p)\in {\mathcal P}^\ast$  if and only if  $r=1$, $p_1\equiv 7\,(\mod{8})$ and $p\equiv g^2\,(\mod{N})$, where $g$ is a generator of $(\Z/N\Z)^\times$. 
\end{proposition}
\proof
It is clear that $r=1$ since  $f$ is odd and $\phi(N)/f=2$. Then, by Theorem~\ref{thm:aoki},  $(N,f,p)\in {\mathcal P}^\ast$ if and only if  $\chi_1^{f_1}(-1)=-1$ and $\chi_1^{f_1}(2)=1$, 
where $\chi_1^{f_1}$ is of order $2$. 
Note that $\chi_1^{f_1}(-1)=-1$ if and only if $p_1\equiv 3\,(\mod{4})$. 
On the other hand, by the supplementary law of quadratic reciprocity, $\chi_1^{f_1}(2)=1$ if and only if $p_1\equiv 1,7\,(\mod{8})$.  

Furthermore, we need to choose $p$ so that $\phi(N)=2\ord_{N}(p)$, i.e., $p\equiv g^2\,(\mod{N})$. 
\qed
\vspace{0.3cm} 

The claim above is also obtainable from the complete characterization of index $2$ Gauss sums~\cite{YX10}. 

\begin{proposition}\label{prop:66}
Assume $\phi(N)/f=6$. Then, $(N,f,p)\in {\mathcal P}^\ast$  if and only if  $r=1$, $p_1\equiv 7\,(\mod{24})$ such that $p_1=a^2+27b^2$ for some integers $a,b$,   and  
$p\equiv g^6\,(\mod{N})$, where $g$ is a generator of $(\Z/N\Z)^\times$. 
\end{proposition}
\proof
Since $f$ is odd and $\phi(N)/2f$ is odd, we have $r=1$. Then, by Theorem~\ref{thm:aoki},  $(N,f,p)\in {\mathcal P}^\ast$ if and only if $\chi_1^{f_1}(-1)=-1$ and $\chi_1^{f_1}(2)=1$, where $\chi_1^{f_1}$ is of order $6$. If $6\not |p_1-1$, it must be $p_1=3$ since $6\,|\,\phi(m_1)$. In this case, $2$ is not a $6$th power modulo $m_1$ since $2$ is a generator of $(\Z/m_1\Z)^\times$. Hence, we have $6|p_1-1$. Note that  $\chi_1^{f_1}(-1)=-1$ if and only if $p_1\equiv 7\,(\mod{12})$. On the other hand, 
by the supplementary law of quadratic reciprocity and the cubic reciprocity law~\cite[Corollary~2.6.4]{BEW97}, $\chi_1^{f_1}(2)=1$ if and only if $p_1\equiv 1,7\,(\mod{8})$ and $p_1=a^2+27b^2$ for some integers $a$ and $b$. 

Furthermore, 
we need to choose $p$ so that $\phi(N)=6\ord_{N}(p)$, i.e., $p\equiv g^6\,(\mod{N})$. 
\qed

\vspace{0.3cm}
The claims (1) and (2) in Theorem~\ref{thm:small} give two sufficient conditions for $(N,f,p)\in {\mathcal P}^\ast$ in the case 
where $\phi(m)/f=4$ and $f$ is odd. 
The two cases in Proposition~\ref{prop:44} below correspond to those two conditions. In particular, we prove that the two conditions are also necessary. 
\begin{proposition}\label{prop:44}
Assume that  $\phi(N)/f=4$. 
Then, $(N,f,p)\in {\mathcal P}^\ast$  if and only if  $r=2$, 
$\gcd{(f_1,f_2)}=1$, and either of the following conditions holds: 
\begin{itemize} 
\item[(1)] $p_1,p_2\equiv 7\,(\mod{8})$; 
\item[(2)] $p_1\equiv 7\,(\mod{8})$, $p_2\equiv 3\,(\mod{4})$ and $p_1$ is quadratic modulo $p_2$. 
\end{itemize}
Furthermore, $p$ is chosen so that $p\equiv g_1^2\,(\mod{2m_1})$ and $p\equiv g_2^2\,(\mod{2m_2})$, where $g_1$ and $g_2$ are generators of $(\Z/2m_1\Z)^\times$ and $(\Z/2m_2\Z)^\times$, respectively. 
\end{proposition}
\proof
Since $\phi(m_i)/f_i\ge 2$, we have $r=1$ or $2$. In the $r=1$ case, by Theorem~\ref{thm:aoki}, $(N,f,p)\in {\mathcal P}^\ast$ if and only if $\chi_1^{f_1}(-1)=-1$ and $\chi_1^{f_1}(2)=1$, where $\chi_1^{f_1}$ is a character of order $4$. Then, $\chi_1^{f_1}(-1)=-1$ if and only if $p_1\equiv 5\,(\mod{8})$. On the other hand, $\chi_1^{2f_1}(2)=1$ if and only if $p_1\equiv 1,7\,(\mod{8})$ by the supplementary law of quadratic reciprocity. Hence, this case is impossible. 

Next, we assume that $r=2$. It is clear that  $\gcd{(f_1,f_2)}=1$; otherwise, $\phi(N)/f>4$. Then, we have 
$\phi(m_i)/f_i=2$, $i=1,2$. In this case, all characters in $\D^-(N,p)$ are $\chi_1^{f_1}$ 
and $\chi_2^{f_2}$, both of which are of order $2$. 
Then, by Theorem~\ref{thm:aoki}, $(N,f,p)\in {\mathcal P}^\ast$  if and only if 
$\chi_1^{f_1}(-1)=\chi_2^{f_2}(-1)=-1$ and either of the following holds: 
$\chi_1^{f_1}(2)=\chi_2^{f_2}(2)=1$, 
$\chi_1^{f_1}(2)=\chi_2^{f_2}(p_1)=1$ (or switching $p_1$ and $p_2$,  $\chi_1^{f_1}(p_2)=\chi_2^{f_2}(2)=1$), or $\chi_1^{f_1}(p_2)=\chi_2^{f_2}(p_1)=1$. 
It is clear that $\chi_1^{f_1}(-1)=\chi_2^{f_2}(-1)=-1$ if and only if  $p_1,p_2\equiv 3\,(\mod{4})$. On the other hand, since $\chi_1^{f_1}(p_2)\chi_2^{f_2}(p_1)=-1$ by the quadratic reciprocity law, the condition that $\chi_1^{f_1}(p_2)=\chi_2^{f_2}(p_1)=1$ is impossible. 
Noting that $\chi_i^{f_i}(2)=1$ if and only if $p_i\equiv 1,7\,(\mod{8})$,  the former two conditions are corresponding to the cases (1) and (2) in the statement, respectively. 
In these cases, noting that $\gcd{(f_1,f_2)}=1$, 
we need to choose $p$ so that  $p \equiv g_1^2\,(\mod{2m_1})$ and $p\equiv g_2^2\,(\mod{2m_2})$. 
\qed

\vspace{0.3cm}
In the following proposition, we treat pure Gauss sums with  $\phi(N)/f= 8$ and $f$ odd, which have not been characterized in the literature. 
\begin{proposition}\label{prop:88}
Assume that  $\phi(N)/f=8$. 
Then,
$(N,f,p)\in \cP^\ast$ if and only if either of the following holds:  
\begin{itemize}
\item[(1)] $r=1$,  $p_1=a^2+64b^2$ for some odd integers $a,b$, and $p\equiv  g^8\,(\mod{N})$, where $g$ is a generator of $(\Z/N\Z)^\times$;  
\item[(2)] $r=2$, $p_1\equiv 5\,(\mod{8})$, $p_2\equiv 3\,(\mod{8})$,   $\gcd{(f_1,f_2)}=1$,  $p_1$ is quadratic modulo $p_2$ and  $p_2$ is quartic modulo $p_1$. Furthermore, $p\equiv g_1^4\,(\mod{2m_1})$ and $p\equiv  g_2^2\,(\mod{2m_2})$, where $g_1$ and $g_2$ are generators of $(\Z/2m_1\Z)^\times$ and $(\Z/2m_2\Z)^\times$, respectively. 
\item[(3)] $r=3$,  $\gcd{(f_i,f_j)}=1$ for any distinct $i,j\in \{1,2,3\}$ and  either one of the following holds: 
\begin{itemize}
\item[i)] $p_1,p_2,p_3\equiv 7 \,(\mod{8})$;  
\item[ii)] $p_1\equiv 7\,(\mod{8})$, $p_2,p_3\equiv 3 \,(\mod{8})$ and  $p_1$ is quadratic modulo  $p_i$ for both $i=2,3$;
\item[iii)] $p_1\equiv 7\,(\mod{8})$, $p_2,p_3\equiv 3 \,(\mod{8})$, $p_1$ is quadratic modulo $p_2$ and $p_2$ is quadratic modulo $p_3$. 
\end{itemize}
Furthermore, $p\equiv g_i^2\,(\mod{2m_i})$ for all $i=1,2,3$,  where $g_i$ is a  generator of $(\Z/2m_i\Z)^\times$, respectively.  
\end{itemize}
\end{proposition}
\proof 
Since $\phi(m_i)/f_i\ge 2$, we have $r=1,2$ or $3$. In the case where $r=1$, by Theorem~\ref{thm:aoki}, $(N,f,p)\in {\mathcal P}^\ast$ if and only if $\chi_1^{f_1}(-1)=-1$ and $\chi_1^{f_1}(2)=1$, where $\chi_1^{f_1}$ is a character of order $8$. Then, $\chi_1^{f_1}(-1)=-1$ if and only if $p_1\equiv 9\,(\mod{16})$. On the other hand, by \cite[Corollary~7.5.8]{BEW97},  $\chi_1^{f_1}(2)=1$ under the assumption that $p_1\equiv 9\,(\mod{16})$ if and only if $p_1=a^2+64b^2$ for some odd integers $a,b$. 
In this case, $p$ must be chosen so that $\phi(N)=8\ord_N(p)$, i.e., $p\equiv g^8\,(\mod{N})$. 


Assume that $r=2$. If $\phi(m_i)/f_i=2$ for $i=1,2$, it follows that $\gcd{(f_1,f_2)}=2$, which contradicts to that $f$ is odd. Hence, we can assume that $\phi(m_1)/f_1=4$, $\phi(m_2)/f_2=2$ and $\gcd{(f_1,f_2)}=1$. Then, all characters in $\D^-(N,p)$  are given as 
\[
\chi_1^{f_1}, \chi_1^{3f_1},  \chi_2^{f_2}, \chi_1^{2f_1}\chi_2^{f_2}.
\]
Note that $\chi_1^{f_1}(-1)=-1$ if and only if $p_1\equiv 5\,(\mod{8})$. Then, by the supplementary law of quadratic reciprocity, we have $\chi_1^{2f_1}(2)=-1$. Since $\chi_1^{2f_1}\chi_2^{f_2}(2)=1$ by Theorem~\ref{thm:aoki}, we have $\chi_2^{f_2}(2)=-1$. Noting that  $\chi_2^{f_2}(-1)=-1$ if and only if $p_2\equiv 3\,(\mod{4})$, by the supplementary law of quadratic reciprocity, we have $p_2\equiv 3\,(\mod{8})$. On the other hand, by 
Theorem~\ref{thm:aoki}, either $\chi_2^{f_2}(2)=1$ or $\chi_2^{f_2}(p_1)=1$ holds. Since $\chi_2^{f_2}(2)=-1$, we have $\chi_2^{f_2}(p_1)=1$. Similarly, we have $\chi_1^{f_1}(p_2)=1$. These conditions correspond to the case (2) in the statement. In this case, $p$ must be chosen so that $p\equiv g_1^4\,(\mod{2m_1})$ and $p\equiv g_2^2\,(\mod{2m_2})$. 

We finally assume that $r=3$. Then, we have $\phi(m_i)/f_i=2$ for every $i=1,2,3$ and $\gcd{(f_i,f_j)}=1$ for any distinct $i,j\in \{1,2,3\}$. In this case, all characters in $\D^-(N,p)$ are 
\[
\chi_1^{f_1},\chi_2^{f_2},\chi_3^{f_3}  \mbox{ and } \chi_1^{f_1}\chi_2^{f_2}\chi_3^{f_3}. 
\] 
Note that $\chi_1^{f_1}(-1)=\chi_2^{f_2}(-1)=\chi_3^{f_3}(-1)=-1$ if and only if $p_1,p_2,p_3\equiv 3\,(\mod{4})$. It follows that  
$\chi_1^{f_1}\chi_2^{f_2}\chi_3^{f_3}(2)=1$ by Theorem~\ref{thm:aoki}. Then, 
by the supplementary law of quadratic reciprocity, $\chi_i^{f_i}(2)=1$ for all $i$ if and only if $p_i\equiv 7\,(\mod{8})$ for all $i$. This corresponds to the case (3)-i) in the statement.  In  other cases, we can assume that  $\chi_1^{f_1}(2)=1$ and  
$\chi_2^{f_2}(2)=\chi_3^{f_3}(2)=-1$. These are equivalent to that $p_1\equiv 7\,(\mod{8})$ and 
$p_2,p_3\equiv 3\,(\mod{8})$, respectively. Furthermore, by Theorem~\ref{thm:aoki}, we have either 
$\chi_2^{f_2}(p_1)=\chi_3^{f_3}(p_1)=1$, 
$\chi_2^{f_2}(p_1)=\chi_3^{f_3}(p_2)=1$ or 
$\chi_2^{f_2}(p_3)=\chi_3^{f_3}(p_2)=1$. Since $\chi_2^{f_2}(p_3)\chi_3^{f_3}(p_2)=-1$ by the quadratic reciprocity law,  $\chi_2^{f_2}(p_3)=\chi_3^{f_3}(p_2)=1$ is impossible. The remaining two conditions are corresponding to (3)-ii) and (3)-iii) in the statement. In these cases,  $p$ must be chosen so that $p\equiv g_i^2\,(\mod{2m_i})$ for $i=1,2,3$.  This completes the proof. 
\qed
\vspace{0.3cm}

We list all $(\overline{p},N)\in {\mathcal P}^\ast_f$ for $N\le 5000$ and odd  $f$ in Tables~\ref{Tab0} and \ref{Tab01} in the appendix. From the computational results, we have the following remark. 
\begin{remark}
For $N\le 5000$ and odd $f$, $(N,f,p)\in {\mathcal P}^\ast$ is in ${\mathcal P}^{(2)}$ or satisfies $\phi(N)/f\le 8$ except for 
$(N,f,\overline{p})=(4042,161,21)$. This exception 
will be characterized in Theorem~\ref{thm:existmany} (see Remark~\ref{rem:except}).
\end{remark}

\section{An application of pure Gauss sums to skew Hadamard difference sets}
\label{sec:appli}
We begin with the following general construction of  skew Hadamard difference sets based on pure Gauss sums. 
\begin{proposition}\label{prop:genecon}
Let $p_i$, $i=1,2,\ldots,r$, be distinct odd primes and $u_i$, $i=1,2,\ldots,r$, be 
positive integers. 
Let $N=2m=2m_1m_2\cdots m_r$, where $m_i=p_i^{u_i}$, $i=1,2,\ldots,r$, and let $p$ be a prime such that $p\equiv 3\,(\mod{4})$. Assume that $(N,f,p)\in {\mathcal P}^\ast$ with $f$ odd.   
Define
\[
Y=\{h>1\mid  h \mbox{ is a divisor of $m$ s.t. $(2h,f,p)\not \in {\mathcal P}$}\}, 
\] 
and 
$I$ as an $m$-subset of $\{0,1,\ldots,N-1\}$ satisfying the following conditions: 
\begin{itemize}
\item[(1)] $\{x\,(\mod{m})\mid x \in I\}=\{0,1,\ldots,m-1\}$. 
\item[(2)] $\sum_{x\in I}\zeta_{2h}^{x}=0$ for any $h\in Y$. 
\end{itemize}
Then, for every odd positive integer $s$, 
\begin{equation}\label{def:DDD}
D=\bigcup_{x \in I}C_i^{(N,p^{fs})}
\end{equation}
forms a skew Hadamard difference set in $(\F_{p^{fs}},+)$. 
\end{proposition}
\proof
First, note that $m\,|\,(p^f-1)/(p-1)$ by Lemma~\ref{lem:evcha}.  Then, $D$ is invariant under the multiplication of $x^2$ for any $x\in \F_p^\ast$. 

Let $\gamma$ be a primitive element of $\F_{p^{fs}}$ and  
let ${\eta'}_{N}$ be a fixed multiplicative character of 
order $N$ of $\F_{p^{fs}}$. 
Furthermore,  let $X$ be the set of all divisors of $m$ and 
$Z$ be the set of odd $1\le j\le N-1$ such that   $N/2\gcd{(j,N)}\in X\setminus Y$. 
Then, by the orthogonality of characters and the conditions (1) and (2), we have for any $a=0,1,\ldots,p^{fs}-2$, 
\begin{align}
\psi_{\F_{p^{fs}}}(\gamma^aD)=&\,\frac{1}{N}\sum_{j=0}^{N-1}\sum_{i\in I}G_{p^{fs}}({\eta'}_N^{j})
{\eta'}_N^{-j}(\gamma^{a+i})\nonumber\\
=&\,\frac{|I|}{N}G_{p^{fs}}({\eta'}_N^{0})+\frac{1}{N}\sum_{j\in Z}\sum_{i\in I}G_{p^{fs}}({\eta'}_N^{j}){\eta'}_N^{-j}(\gamma^{a+i}). \label{eq:totyu}
\end{align} 
Since $G_{p^{fs}}({\eta'}_N^{j})$ is pure for any $j\in Z$, we have 
\[
\sum_{j\in Z}\sum_{i\in I}G_{p^{fs}}({\eta'}_N^{j}){\eta'}_N^{-j}(\gamma^{a+i})\equiv 0\,(\mod{p^{\frac{fs-1}{2}}}). 
\]
Finally, noting that  $G_{p^{fs}}({\eta'}_N^{0})=-1$ and $\gcd{(m,p)}=1$, we have $2\psi_{\F_{p^{fs}}}(\gamma^aD)+1\equiv 0\,(\mod{p^{\frac{fs-1}{2}}})$. Then, by Lemma~\ref{lem:skewcheck}, the claim follows. \qed
\vspace{0.3cm}

Even if we determine the set $Y$, i.e., for which $h$ we have 
$(2h,f,p)\not\in {\mathcal P}$, in the proposition above, 
it may happen that there is no nontrivial subset $I$ satisfying the conditions (1) 
and (2) as commented in Section~\ref{sec:conc}. 
Moreover, it is difficult to determine the dual of $D$ in general. Indeed, to do this, we need to evaluate the signs (or roots of unity) of the corresponding pure Gauss sums. Thus, we have to choose suitable $(N,f,p)\in {\mathcal P}^\ast$ such that a nontrivial subset $I$ exists satisfying the conditions (1) and (2) and we can evaluate the signs (or roots of unity)  in some sense. 
From this point of view, we consider pure Gauss sums satisfying a special property defined below. 

Throughout this section, we assume that the order $f$ of $p$ modulo $N$ is odd and $N$ has the prime factorization $N=2m_1m_2\cdots m_r$, where $m_i=p_i^{u_i}$ with $p_i$ an odd prime and $u_i$ a positive integer. 
We consider $(N,f,p)\in {\mathcal P}^\ast$ such that 
\begin{itemize}
\item[($\star$)] $(2\prod_{i\in J}m_i, f, p)\in {\mathcal P}$ for any subset $J\subseteq \{1,2,\ldots,r\}$ such that $1\in J$. 
\end{itemize}
\subsection{Characterization of pure Gauss sums with property ($\star$)}
In this subsection, we give a characterization of pure Gauss sums with property ($\star$). 
\begin{lemma}\label{lem:red1}
If $(N,f,p)\in {\mathcal P}^\ast$ and $(2\prod_{i\in J}m_i,f,p)\in {\mathcal P}$ 
for a subset $J\subseteq \{1,2,\ldots,r\}$, then $(2\prod_{i\in J}m_i,f',p)\in {\mathcal P}^\ast$, where $f'=\lcm{(f_i: i \in J)}$. 
\end{lemma}
\proof
Note that $f=\lcm{(f_1,f_2,\ldots,f_r)}$ and $f'\,|\,f$. Then, 
by Theorem~\ref{thm:lift}, we have $(2\prod_{i\in J}m_i,f',p)\in {\mathcal P}^\ast$. 
\qed
\begin{proposition}\label{prop:nece_am}
Assume that $(N,f,p)\in {\mathcal P}^\ast$ with property {\em ($\star$)}. 
Then, either of the following holds: 
\begin{itemize}
\item[(1)] $\chi_i^{f_i}(2)=1$ for any $i\in \{1,2,\ldots,r\}$; or  
\item[(2)] $\chi_1^{f_1}(2)=1$, $\chi_i^{f_i}(2)\not=1$ for any $i\in \{2,3,\ldots,r\}$, and $\chi_i^{f_i}(p_1)=1$ for any $i\in \{2,3,\ldots,r\}$. 
\end{itemize}
\end{proposition}
\proof 
Since $(2m_1,f,p)\in \cP$,  we have $(2m_1,f_1,p)\in \cP^\ast$ by Lemma~\ref{lem:red1}. Then, by Theorem~\ref{thm:aoki}, we have 
$\chi_{1}^{f_1}(2)=1$. Furthermore, 
since $(2m_1m_i,f,p)\in \cP$ for any $i\in \{2,3,\ldots,r\}$, we have 
$(2m_1m_i,\lcm{(f_1,f_i)},p)\in \cP^\ast$ by Lemma~\ref{lem:red1} again. Then, 
by Theorem~\ref{thm:aoki}, we have either $\chi_i^{f_i}(2)=1$ or $\chi_i^{f_i}(p_1)=1$. 

We assume that $\chi_i^{f_i}(2)=1$ for some $i\in \{2,3,\ldots,r\}$. 
Since $(2m_1m_im_j,f,p)\in \cP$ for any $j\in \{2,3,\ldots,r\}\setminus\{i\}$, 
$(2m_1m_im_j,\lcm{(f_1,f_i,f_j)},p)\in \cP^\ast$ follows by Lemma~\ref{lem:red1}. Then, by Theorem~\ref{thm:aoki}, we have 
$\chi_{1}^{f_1}\chi_{i}^{f_i}\chi_{j}^{f_j}(2) =1$. Since $\chi_{1}^{f_1}(2)=\chi_{i}^{f_i}(2)=1$, we have $\chi_{j}^{f_j}(2)=1$. This implies that $\chi_{h}^{f_h}(2)=1$ for any $h\in \{1,2,\ldots,r\}$. If $\chi_i^{f_i}(2)\not=1$ for any $i\in \{2,3,\ldots,r\}$, we have $\chi_i^{f_i}(p_1)=1$ for any $i\in \{2,3,\ldots,r\}$. This completes the proof. 
\qed
\begin{remark}
If $(N,f,p)$ satisfies the condition (2) in Proposition~\ref{prop:nece_am},  
$(2m_1,f,p) \in {\mathcal P}$ but 
$(2m_i,f,p)\not \in {\mathcal P}$ for any $i\in \{2,3,\ldots,r\}$ by Theorem~\ref{thm:aoki}. 
\end{remark}
We now give a sufficient condition for $(N,f,p)\in \cP^\ast$ with property ($\star$). 
\begin{theorem}\label{thm:existmany1}
Assume that $f=\ord_{N}(p)$ is odd. 
If $(N,f,p)\in {\mathcal P}^{(2)}$,  $(N,f,p)\in \cP^\ast$ with property~{\em ($\star$)}. 
\end{theorem}
\proof
Since $2\in \langle p\rangle\,(\mod{m'})$ for any divisor $m'$ of $m$, the assertion holds.  
\qed
\vspace{0.3cm}

The theorem above implies that  $(N,f,p)\in {\mathcal P}^{(2)}$ belongs to 
the class (1) of Proposition~\ref{prop:nece_am}. The pure Gauss sums in this case were used for constructing skew Hadamard difference sets as in Theorem~\ref{thm:fx}. Next, we give a sufficient condition for $(N,f,p)\in \cP^\ast$ to belong to the 
 class (2) of Proposition~\ref{prop:nece_am}.
\begin{theorem}\label{thm:existmany}
Assume that $f_i$ are all odd and $\gcd{(f_1,f_i)}=1$ for any $i\in \{2,3,\ldots,r\}$.  
If $(N,f,p)$ satisfies that $\phi(m_1)/f_1=2$, $2\in \langle p\rangle\,(\mod{m_1})$, $-2\in \langle p\rangle\,(\mod{m/m_1})$ and  
$p_1\in \langle p\rangle\,(\mod{m/m_1})$, 
then $(N,f,p)\in \cP^\ast$ with  property~{\em ($\star$)}. 
\end{theorem}
\proof 
By the assumption that $\gcd{(f_1,f_i)}=1$ for any $i\in \{2,3,\ldots,r\}$, any odd character in $\D^-(N,p)$ has the form 
$\tau_1=\chi_1^{f_1}\chi$ for some even character $\chi$ modulo $m/m_1$ such that $\chi(p)=1$ or $\tau_2=\chi'$ for some odd character $\chi'$ modulo $m/m_1$ such that $\chi'(p)=1$.  
By the assumptions that $\phi(m_1)/f_1=2$, $2\in \langle p\rangle\,(\mod{m_1})$ and $-2\in \langle p\rangle\,(\mod{m/m_1})$,  we have 
\[
\tau_1(2)=\chi_1^{f_1}(2)\chi(2)=\chi(-1)\chi(p^i)=1
\]
for some $i$. On the other hand, since $p_1\in \langle p\rangle\,(\mod{m/m_1})$, we have 
\[
\tau_2(p_1)=\chi'(p_1)=1. 
\]
%
Then, by Theorem~\ref{thm:aoki}, it follows that $(N,f,p)\in \cP^\ast$. Furthermore, it is clear that the property~($\star$) is satisfied. 
\qed
\vspace{0.3cm}

The theorem above is a generalization of Proposition~\ref{prop:44}~(2). 
\begin{remark}\label{rem:except}
The  exception $(N,\overline{p})=(4042,21)\in {\mathcal P}_{161}^\ast$ listed in Table~\ref{Tab01} satisfies the condition in Theorem~\ref{thm:existmany} as 
$m_1=47$ and $m_2=43$.   
\end{remark}

\begin{corollary}\label{cor:coco}
Assume that  
\begin{itemize}
\item[(1)] $p_1\equiv 7\,(\mod{8})$ and $p_i\equiv 3\,(\mod{8})$ for  all $i\in \{2,3,\ldots,r\}$; 
\item[(2)] $p_1$ is quadratic modulo $p_i$ for all $i\in \{2,3,\ldots,r\}$;  
\item[(3)] $f_i=\phi(m_i)/2$ for all $i$; 
\item[(4)] $\phi(m_i)/2$'s are mutually coprime. 
\end{itemize}
Then, $(N,f,p)\in \cP^\ast$ with property {\em ($\star$)}. 
\end{corollary}
\proof
By the supplementary law of  quadratic reciprocity, we have 
$\chi_1^{f_1}(2)=1$ and $\chi_i^{f_i}(2)=-1$ for all $i>1$. 
Then, the conditions of Theorem~\ref{thm:existmany} are fulfilled.  
\qed

\vspace{0.3cm}
One can see that there are infinitely many tuples of $m_1,m_2,\ldots,m_r,p$ 
satisfying the condition of Corollary~\ref{cor:coco}.  
\begin{example}
Fix $m_2=3$ and $m_3=11$. Let $m_1\equiv 7\,(\mod{8})$ be a prime such that $m_1$ is  quadratic modulo both $3$ and $11$, i.e., $m_1\equiv 1\,(\mod{3})$ and 
$m_1\equiv 1,3,4,5,9\,(\mod{11})$. Furthermore, we need the restriction $m_1\equiv 3,5,7,9\,(\mod{10})$ in order to satisfy that $\gcd{(\phi(m_1)/2,\phi(11)/2)}=1$. There are infinitely many such primes $m_1$ by the Dirichlet theorem in arithmetic progressions. For example, we can take $m_1=103,199,223,367,463,\ldots$. Let $p'$ be any integer with $1\le p'\le N-1$ determined by the congruences 
\[
\begin{cases}
p'\equiv 1\,(\mod{6})\\
p'\equiv 1,3,5,9 \mbox{ or }15\,(\mod{22}) \\
p'\equiv g^2\,(\mod{2m_1})
\end{cases}
\]
for any generator $g$ of $(\Z/2m_1\Z)^\times$. Then, $\langle p'\rangle$ is of index $2$ modulo  $2m_i$ for each $i$. Then, $(N,f,p)\in \cP^\ast$ with property {\em ($\star$)} for any odd prime $p\equiv p'\,(\mod{N})$, where $f:=5(m_1-1)/2$.   
\end{example}

\subsection{The signs  or roots of unity of pure Gauss sums with property ($\star$)}
 In this subsection, we study the signs or roots of unity of pure Gauss sums satisfying property ($\star$).
The following result was known. 
\begin{lemma}\label{lem:signAA}{\em (\cite[Lemma~6]{E81})}
If $G_q(\eta_N)$ is pure, then $\epsilon=G_q(\eta_N)/q^{\frac{1}{2}}$ is a $2\gcd{(N,p-1)}$th root of unity. 
\end{lemma} 
The sign (root of unity) ambiguities of pure Gauss sums in the class ${\mathcal P}^{(2)}$ was completely 
determined  as in Theorem~\ref{prop:2m}. 
However,  it is difficult to explicitly determine them for pure Gauss sums in general. In fact, it sometimes becomes complicated as in 
\cite[Theorem~10]{E81}. In this subsection,  we show that 
if $G_{p^f}(\eta_N)$ is pure with property $(\star)$, 
the sign (or root of unity) of  
$G_{p^f}(\eta_N)$  is determined from those of $G_{p^f}(\eta_{2m_1m_i})$'s, $1\le i\le r$. This property will be used to determine the duals of skew Hadamard difference sets obtained from the construction in Theorem~\ref{thm:mainskew}.  
%

For positive integers $x$ and $y$ with $\gcd{(x,y)}=1$, let $\inv(x;y)$ denote an integer such that $x\cdot \inv(x;y)\equiv 1\,(\mod{y})$. The claim of the following lemma was given  in the proof of  \cite[Theorem~7]{E81}. 
\begin{lemma}\label{lem:sign}
Let $n=st$ be a positive integer, where $s$ is an odd prime power and $t>1$ with $\gcd{(s,t)}=1$. Let $\eta_n$ be a multiplicative character 
of order $n$ of $\F_{p^f}$. 
Assume that $(n,f,p)\in \cP$ and $(t,f,p)\in \cP$. Then, $G_{p^f}(\eta_{n})/G_{p^f}(\eta_{n}^{s\cdot \inv(s,t)})$ is a $\gcd{(p-1,s)}$th root of unity. 
\end{lemma}
Recall that $N=2m=2m_1m_2\cdots m_r$, where $m_i=p_i^{u_i}$, $i=1,2,\ldots,r$. For a subset $J\subseteq \{1,2,\ldots,r\}$, denote $m_J=\prod_{i \in J}m_i$ and $n_J=m/m_J$. Let $\omega$ be a primitive root of $\F_{p^f}$, and let $\eta_{N}$ be a fixed multiplicative character of order $N$ of $\F_{p^f}$ such that $\eta_{N}(\omega)=\zeta_{N}$. We denote  $\eta_{2m_J}=\eta_{N}^{n_J\cdot\inv(n_J,2m_J)}$, and also denote an arbitrary multiplicative character of order $h$ of $\F_{p^f}$ by 
$\theta_h$. Note that for any $m_i\,|\,m/m_J$, 
$\eta_{2m_J}=\eta_{2m_im_J}^{m_i\cdot \inv(m_i,2m_J)}$. 
\begin{proposition}\label{prop:sign_c1}
Assume that $(N,f,p)\in \cP^\ast$ with property~{\em ($\star$)}. 
Then, the following hold: 
\begin{itemize}
\item  There are integers $s_i$, $i=1,2,\ldots,r$, such that
\begin{equation}\label{eq:sign21}
G_{p^f}(\eta_{2m_1})=\zeta_{m_1}^{s_1}G_{p^f}(\eta_2)\mbox{\, and \, }
G_{p^f}(\eta_{2m_1m_i})=\zeta_{m_i}^{s_i}G_{p^f}(\eta_{2m_1}), \quad 2\le i\le r.  
\end{equation} 
\item 
Let $J\subseteq \{1,2,3,\ldots,r\}$ such that $1\in J$. Then, 
\begin{equation}\label{eq:sign3}
G_{p^f}(\eta_{2m_J})=\Big(\prod_{i\in J}\zeta_{m_i}^{s_i}\Big) G_{p^f}(\eta_2), 
\end{equation}
where $s_i$'s are defined as in \eqref{eq:sign21}. In particular, $s_1=0$. 
\end{itemize}
\end{proposition}
\proof 
\eqref{eq:sign21} is a direct consequence of Lemma~\ref{lem:sign}. 
Furthermore, we have $s_1=0$ by Proposition~\ref{prop:2m}. 
 
We prove that \eqref{eq:sign3} holds. 
Let $m_i,m_j$ be distinct prime power divisors of  $m/m_1$ 
and  $J$ be any subset of $\{1,2,\ldots,r\}\setminus \{i,j\}$ containing $1$.  
By Lemma~\ref{lem:sign}, we have
\[
\begin{cases} 
G_{p^f}(\eta_{2m_im_J})/G_{p^f}(\eta_{2m_J})=\zeta_{m_i}^{u}& \\ 
G_{p^f}(\eta_{2m_jm_J})/G_{p^f}(\eta_{2m_J})=\zeta_{m_j}^{u'}& 
\end{cases}\mbox{ and \, }
\begin{cases} 
G_{p^f}(\eta_{2m_im_jm_J})/G_{p^f}(\eta_{2m_im_J})=\zeta_{m_i}^v& \\ 
G_{p^f}(\eta_{2m_im_jm_J})/G_{p^f}(\eta_{2m_jm_J})=\zeta_{m_j}^{v'}&
\end{cases}
\]
for some integers $u,u',v,v'$.  By combining these equations, we have 
$\zeta_{m_i}^{u}\zeta_{m_j}^{v'}=\zeta_{m_i}^v\zeta_{m_j}^{u'}$, i.e., $u=v$ and 
$u'=v'$. Hence, $G_{p^f}(\eta_{2m_im_jm_J})/G_{p^f}(\eta_{2m_J})=\zeta_{m_i}^{u}\zeta_{m_j}^{u'}$. This argument inductively shows \eqref{eq:sign3}. 
\qed


\begin{proposition}\label{prop:sign_c2}
Assume that $(N,f,p)\in \cP^\ast$ with property~{\em ($\star$)}. 
Let $J\subseteq \{1,2,\ldots,r\}$ such that $1\in J$, and let $J'\subseteq J$. 
Furthermore, let $v=\prod_{i\in J'}p_i^{\ell_i}$ for some integers $1\le \ell_i\le u_i-1$.
Then,  
\begin{equation}\label{eq:sign31}
G_{p^f}(\eta_{2m_J}^v)=\Big(\prod_{i\in J}\zeta_{m_i}^{vs_i}\Big) G_{p^f}(\eta_2). 
\end{equation}
\end{proposition}
\proof 
By Corollary~\ref{cor:HD}, 
we have 
\begin{equation}\label{eq:hasse1}
G_{p^f}(\eta_{2m_J}^v)=
p^{-f\frac{v-1}{2}}\sigma_{1,v^{-1}}\Big(\prod_{j=0}^{v-1}G_{p^f}(\eta_{2m_J}\theta_v^j)\Big),  
\end{equation}
where $\sigma_{1,v^{-1}}\in \Gal(\Q(\zeta_{2m_J},\zeta_{p})/\Q)$ 
and $\theta_v$ is any multiplicative character of order $v$ of $\F_{p^f}$.  Noting that $\gcd{(2m_J,1+2m_Jj/v)}=1$ for any $j=0,1,\ldots,v-1$, we have  $\sigma_{1+2m_Jj/v,1}(G_{p^f}(\eta_{2m_J}))=G_{p^f}(\eta_{2m_J}^{1+2m_Jj/v})$. On the other hand, by Proposition~\ref{prop:sign_c1}, we have $\sigma_{1+2m_Jj/v,1}(G_{p^f}(\eta_{2m_J}))=\Big(\prod_{i\in J}\zeta_{m_i}^{s_i(1+2m_J j/v)}\Big) G_{p^f}(\eta_2)$. Hence, $G_{p^f}(\eta_{2m_J}^{1+2m_Jj/v})=\Big(\prod_{i\in J}\zeta_{m_i}^{s_i(1+2m_J j/v)}\Big) G_{p^f}(\eta_2)$. 
Then, by noting that $v$ is odd, 
\begin{align*}
\prod_{j=0}^{v-1}
G_{p^f}(\eta_{2m_J}\theta_v^j)=&\, \prod_{j=0}^{v-1}G_{p^f}(\eta_{2m_J}^{1+2m_Jj/v})\\
=&\, \Big(\prod_{i\in J}\zeta_{m_i}^{s_iv}\Big)\Big(\prod_{i\in J}\zeta_{m_i}^{s_i2m_J (1+\cdots+v-1)/v} \Big)
G_{p^f}(\eta_2)^v\\
=&\,\Big(\prod_{i\in J}\zeta_{m_i}^{s_iv} \Big)
G_{p^f}(\eta_2)^v. 
\end{align*}
Furthermore, we have 
\[
G_{p^f}(\eta_2)^v=\eta_2^{\frac{v-1}{2}}(-1)p^{f\frac{v-1}{2}}G_{p^f}(\eta_2). 
\]
%
Therefore, 
\eqref{eq:hasse1} is reformulated as 
\[
G_{p^f}(\eta_{2m_J}^{v})=
\eta_{2}^{\frac{v-1}{2}}(-1) \Big(\prod_{i\in J}\zeta_{m_i}^{s_iv} \Big) \sigma_{1,v^{-1}}(G_{p^f}(\eta_{2}))
=\eta_{2}^{\frac{v-1}{2}}(-1)\eta_2(v) \Big(\prod_{i\in J}\zeta_{m_i}^{s_iv} \Big) G_{p^f}(\eta_{2}). 
\]
Finally, we see that $\eta_{2}^{\frac{v-1}{2}}(-1)
\eta_{2}(v)=1$.  
Note that $\eta_2(v)=\prod_{i\in J'}\big(\frac{p_i}{p}\big)^{\ell_i}$, where $\big(\frac{p_i}{p}\big)$ is the Legendre symbol. Since $f$ is odd, we have 
$\big(\frac{p}{p_i}\big)=1$. 
Then, by the quadratic reciprocity law, we have 
\[
\prod_{i\in J'}\Big(\frac{p_i}{p}\Big)^{\ell_i}=\prod_{i\in J'}(-1)^{\frac{(p-1)(p_i-1)\ell_i}{4}}. 
\]
Let $h$ be the number of $i\in J'$ such that $p_i\equiv 3\,(\mod{4})$ and $\ell_i$ is odd. 
Then, we have $\eta_2(v)=(-1)^{h(p-1)/2}$. On the other hand, 
$\eta_{2}^{\frac{v-1}{2}}(-1)=1$ if and only if $p\equiv 1\,(\mod{4})$ or $p\equiv 3\,(\mod{4})$ and $v\equiv 1\,(\mod{4})$ (i.e., $h$ is even). Hence,  $\eta_2(v)=\eta_{2}^{\frac{v-1}{2}}(-1)$. 
This completes the proof of the proposition. 
\qed

\begin{remark}\label{rem:signs1}
\begin{itemize}
\item[(1)] 
Recall that $\eta_{2m_J}=\eta_{2m}^{n_J\cdot \inv(n_J,2m_J)}$. 
Let 
\begin{equation}\label{def:a}
A=s_1m_2\cdots m_{r}+\cdots +s_r m_1\cdots m_{r-1}.
\end{equation} 
Then, we have 
\[
\eta_{2m_J}(\omega^{2A})=\prod_{i=1}^r \zeta_{m_i}^{s_i n_J\cdot \inv(n_J,2m_J)}
=\prod_{i\in J} \zeta_{m_i}^{s_i n_J\cdot \inv(n_J,2m_J)}. 
\]
Since $n_J\cdot \inv(n_J,2m_J)\equiv 1\,(\mod{m_i})$ for $i\in J$, we have 
$\eta_{2m_J}(\omega^{2A})=\prod_{i\in J} \zeta_{m_i}^{s_i}$. Hence, by Propositions~\ref{prop:sign_c1} and \ref{prop:sign_c2}, for any odd $j$ such that 
$p_1\,|\,\frac{N}{\gcd{(j,N)}}$, 
\begin{equation}\label{eq:finaleq}
G_{p^f}(\eta_{N}^j)=\eta_{N}^j(\omega^{2A})G_{p^f}(\eta_2). 
\end{equation}
\item[(2)] If $(N,p,f)$ satisfies the condition of   Theorem~\ref{thm:existmany}, 
$\eta_N(\omega^{2A})$ is a cubic root of unity. In fact, the condition $-2\in \langle p\rangle\,(\mod{m/m_1})$ implies that $\gcd{(m_i,p-1)}=1$ or $3$ for any $i=2,3,\ldots,r$. Then,  by Lemma~\ref{lem:sign} and Proposition~\ref{prop:sign_c1}, the claim follows. 
\end{itemize}
\end{remark}

\subsection{A construction of skew Hadamard difference sets}\label{sec:const}
In this subsection, we show that if $(N,f,p)\in {\mathcal P}^\ast$ with property 
$(\star)$, there are nontrivial choices of $I$ satisfying 
the conditions (1) and (2) of Proposition~\ref{prop:genecon}. Furthermore, we can determine the dual of $D$ in this case. First, we illustrate our construction giving one example below.  
\begin{example}\label{exam:skew}
As in Table~\ref{Tab0}, we have $(42,3,67)\in \cP^\ast$ and $(14,3,67)\in \cP$, i.e., it satisfies the condition~{\em ($\star$)} as $m_1=7$. Note that 
$(6,3,67)\not\in \cP$, i.e., $(42,3,67)\in \cP^\ast$ belongs to the class (2) of Proposition~\ref{prop:nece_am}. Let $I$ be any $21$-subset of $\{0,1,\ldots,41\}$ satisfying the following conditions: 
\begin{itemize}
\item[(1)] $\{x\,(\mod{21})\mid x \in I\}=\{0,1,\ldots,20\}$;  
\item[(2)] $\sum_{x\in I}\zeta_6^{x}=0$. 
\end{itemize}
For example, we can take 
\[
I=\{0, 1, 2, 3, 4, 5, 6, 7, 8, 9, 10, 11, 12, 13, 14, 15, 16, 17, 18, 20, 40\}.
\]
Then, for every  odd positive integer $s$, 
\[
D=\bigcup_{x \in I}C_i^{(42,67^{3s})}
\]
forms a skew Hadamard difference set in $(\F_{67^{3s}},+)$. 
Furthermore, its dual is given as $D^\perp=\{\psi_a\in G^\perp\,|\,a\in \bigcup_{i\in I'}C_i^{(N,p^{fs})}\}$, where $I'=-I+14s$. 
\end{example}

\begin{theorem}\label{thm:mainskew}
With notations as in Proposition~\ref{prop:genecon}, assume that 
$(N,f,p)\in \cP^\ast$ with property~{\em ($\star$)}. Furthermore, redefine $Y$ as 
\[
Y=\{h>1\mid  h \mbox{ is a divisor of $\prod_{i=2}^rm_i$}\}. 
\] 
Then, for every odd positive integer $s$, the set $D$ defined in \eqref{def:DDD} 
forms a skew Hadamard difference set in $(\F_{p^{fs}},+)$. In particular, its dual  is given as $D^\perp=\{ \psi_a\in G^\perp\,|\,a\in \bigcup_{i\in I'}C_i^{(N,p^{fs})}\}$, where $I'=-I+2As$, 
where $A$ is defined as in \eqref{def:a}.
\end{theorem}
\proof
Let $\gamma$ be a primitive element of $\F_{p^{fs}}$ and let  $\omega=\gamma^{(p^{fs}-1)/(p^f-1)}$. Furthermore, let $\eta_{N}$ be a fixed multiplicative character of 
order $N$ of $\F_{p^f}$ such that $\eta_{N}(\omega)=\zeta_{N}$, and let $\eta_N'$ be the lift of $\eta_N$ to $\F_{p^{fs}}$. 
Continuing from \eqref{eq:totyu}, we have by Theorem~\ref{thm:lift} that 
\begin{align}
\psi_{\F_{p^{fs}}}(\gamma^aD)
=-\frac{1}{2}+\frac{1}{N}\sum_{j\in Z}\sum_{i\in I}(G_{p^{f}}(\eta_{N}^j))^s
\eta_{N}^{-j}(\omega^{a+i}).   \label{eq:skew1}
\end{align} 
By Proposition~\ref{prop:sign_c2} and Remark~\ref{rem:signs1}, for any $j\in Z$ 
\[
G_{p^f}(\eta_{N}^j)=\eta_{N}^j(\omega^{2A}) G_{p^f}(\eta_2). 
\] 
Hence, 
continuing from \eqref{eq:skew1}, we have 
\begin{align*}\label{eq:skew2}
\psi_{\F_{p^{fs}}}(\gamma^aD)=&\,-\frac{1}{2}+\frac{(G_{p^f}(\eta_2))^s}{N}\sum_{j \in Z}\sum_{i\in I}
\eta_{N}^{-j}(\omega^{a+i})\eta_{N}^{js}(\omega^{2A})\\
=&\,-\frac{1}{2}+\frac{G_{p^{fs}}(\eta_2')}{N}\sum_{j=1}^{N-1}\sum_{i\in I}
\eta_{N}^{-j}(\omega^{a-2As+i})\\
=&\,-\frac{1}{2}+\frac{G_{p^{fs}}(\eta_2')}{N}\cdot 
\begin{cases}
m& \mbox{ if $-a+2As\in I$,}\\
-m& \mbox{ otherwise,}
\end{cases}
\end{align*}
where $\eta_2'$ is the quadratic character of $\F_{p^{fs}}$.  
Hence, by \eqref{eq:Gaussquad}, we obtain $\psi_{\F_{p^{fs}}}(\gamma^aD)\in  \{\frac{-1\pm \sqrt{-p^{fs}}}{2}\}$. This implies that $D$ is a skew Hadamard difference set. Furthermore, its dual is determined as desired. 
\qed

\begin{remark}
There are nontrivial choices of $I$ satisfying the conditions of Theorem~\ref{thm:mainskew}. For example, let $S_i=\{2im/m_1+j\,|\,0\le j\le 2m/m_1-1\}$ for $i=0,1,\ldots,(m_1-3)/2$ and let   $A_0,A_1$ be an arbitrary partition of $J=\{0,1,\ldots,(m_1-3)/2\}$.  
Then, we can take 
\begin{align*}
I=&\,\left\{x\mid x\in \bigcup_{i\in A_0} S_i\right\}
\cup \left\{x+m\mid x\in \bigcup_{i\in A_1} S_i\right\}\\
& \cup \left\{\frac{(m_1-1)m}{m_1}+2i\mid i=0,1,\ldots,\frac{m/m_1-1}{2}\right\}\\ 
&\cup \left\{\frac{(m_1-1)m}{m_1}+m+2i-1\mid i=1,2,\ldots,\frac{m/m_1-1}{2}\right\}.
\end{align*}
This is a generalization of $I$  in Example~\ref{exam:skew}. 
\end{remark}

The following result is immediately obtained 
by applying  Theorem~\ref{thm:mainskew} to the class ${\mathcal P}^{(2)}$. 
\begin{corollary}\label{cor:const1}
Assume that $2\in \langle p\rangle \,(\mod{m})$. 
Let $Y=\{h>1\mid  h \mbox{ is a divisor of $\prod_{i=2}^rm_i$}\}$, and let $I$ be an arbitrary $m$-subset of $\{0,1,\ldots,N-1\}$ such that $\{x\,(\mod{m})\mid x \in I\}=\{0,1,\ldots,m-1\}$ and   
$\sum_{x\in I}\zeta_{2h}^{x}=0$ for any 
$h\in Y$.  
Then, for any odd positive integer $s$, 
$
D=\bigcup_{x \in I}C_i^{(N,p^{fs})}$
forms a skew Hadamard difference set in $(\F_{p^{fs}},+)$. 
\end{corollary}
\proof
The assumption that $2\in \langle p\rangle \,(\mod{m})$ implies that $(N,f,p)\in {\mathcal P}^{(2)}$ with property ($\star$) belonging to the class (1) of Proposition~\ref{prop:nece_am} by Theorem~\ref{thm:existmany1}. Then, by Theorem~\ref{thm:mainskew}, the claim follows. 
\qed
\vspace{0.3cm}

Note that the result above is contained in Theorem~\ref{thm:fx}. In fact, the construction given  in Theorem~\ref{thm:fx} allows $I$ as an arbitrary subset of $\{0,1,\ldots,N-1\}$ satisfying the condition (1) but not necessarily satisfying the condition (2) since $(2m',f,p)\in {\mathcal P}^\ast$ for any divisor $m'$ of $m$. This also follows from Proposition~\ref{prop:genecon} as $Y=\emptyset$. 

The following  is a new result not within the framework of  Theorem~\ref{thm:fx}.  
\begin{corollary}\label{cor:const2}
Assume that $f_i$ are all odd and $\gcd{(f_1,f_i)}=1$ for any $i\in \{2,3,\ldots,r\}$.  
If $(N,f,p)$ satisfies that $\phi(m_1)/f_1=2$, $2\in \langle p\rangle\,(\mod{m_1})$, $2\in -\langle p\rangle\,(\mod{m/m_1})$ and  
$p_1\in \langle p\rangle\,(\mod{m/m_1})$, 
Let $Y=\{h>1\mid  h \mbox{ is a divisor of $\prod_{i=2}^rm_i$}\}$, and let $I$ be an arbitrary $m$-subset of $\{0,1,\ldots,N-1\}$ such that $\{x\,(\mod{m})\mid x \in I\}=\{0,1,\ldots,m-1\}$ and   
$\sum_{x\in I}\zeta_{2h}^{x}=0$ for any 
$h\in Y$.  
Then, for any odd positive integer $s$, 
$
D=\bigcup_{x \in I}C_i^{(N,p^{fs})}$
forms a skew Hadamard difference set in $(\F_{p^{fs}},+)$. 
\end{corollary}
\proof 
The assumptions imply that $(N,f,p)\in {\mathcal P}^\ast$ with property ($\star$) belonging to the class (2) of Proposition~\ref{prop:nece_am} by 
Theorem~\ref{thm:existmany}.  
Then, by Theorem~\ref{thm:mainskew}, the claim follows.  
\qed

\begin{remark}
We remark that any $(N,f,p)\in {\mathcal P}^\ast$  satisfying the condition of  Proposition~\ref{prop:44}~(2) or  \ref{prop:88}~(3)-ii) 
also 
satisfies the condition of  Corollary~\ref{cor:const2}. Hence, by Tables~\ref{Tab0} and \ref{Tab01}, there exist $(N,f,p)\in {\mathcal P}^\ast$ satisfying the condition of Corollary~\ref{cor:const2}  in abundance.
\end{remark}

\begin{remark}
In this remark, we discuss the inequivalence problem  on  skew Hadamard difference sets obtained from Theorem~\ref{thm:mainskew}. 
Two skew Hadamard difference sets $D_1$ and $D_2$ in an abelian group $G$ are called {\it equivalent} if there exists an automorphism $\sigma\in \Aut(G)$ and an element $x\in G$ such that $\sigma(D_1)+x=D_2$.

Let $D$ be  a skew Hadamard difference set in $(\F_q,+)$.  For a fixed $a\in \F_p^\ast$, define
\[
T_{x,a}(D):=|D\cap (D-x)\cap (D-a\cdot x)|, \quad x\in \F_q^\ast,
\]
and 
\[
n_a(D)=|\{T_{x,a}(D)\mid x\in \F_{q}^\ast\}|. 
\]
It is known that $n_a(D)$ is  an invariant of the equivalence of skew Hadamard difference sets in $(\F_q,+)$, cf.~\cite{Mo1}. 

It is clear that the Paley difference set $D_P$ satisfies that $n_a(D_P)\le 2$ for any $a\in \F_p^\ast$. 
If a skew Hadamard difference set $D$ satisfies $n_a(D)\ge 3$ for some $a\in \F_p^\ast$,  then $D$ is inequivalent to $D_P$. 
Let $D_0\subseteq \F_{67^3}$ be the skew Hadamard difference set in Example~\ref{exam:skew}. We checked by  a computer that 
$n_a(D_0)\ge 3$ for $a=3$.  
On the other hand, since $(N,f,p)=(14,3,67)\in {\mathcal P}^{(2)}$, the set 
$D=\bigcup_{i \in I}C_{i}^{(14,67^3)}$ is also a skew Hadamard difference set for any $7$-subset $I$ of $\{0,1,\ldots,13\}$ such that $I\cap \{x+7\,(\mod{14})\mid x\in I\}=\emptyset$. Let
\begin{align*}
& I_1=\{0,1,2,3,4,5,6\},\,   I_2=\{0, 1, 2, 3, 4, 6,12\},\\
&I_3=\{0, 1, 6, 9, 10, 11, 12\},\, I_4=\{0,1, 2, 4, 6, 10, 12\}
\end{align*}
and define $D_j=\bigcup_{i \in I_j}C_{i}^{(14,67^3)}$ for $j=1,2,3,4$. 
We checked by a computer that $D_j$, $j=1,2,3,4$, are mutually inequivalent 
and they are also inequivalent to the Paley difference set.  Furthermore, it holds that  $n_3(D_0)\not=n_3(D_j)$ for any $j=1,2,3,4$.  
Hence, $D_0$ is inequivalent to $D_j$'s. Thus, 
Corollary~\ref{cor:const2} can give rise to skew Hadamard difference sets 
not obtained from Theorem~\ref{thm:fx}. 
\end{remark}

\section{Concluding remarks}\label{sec:conc}
In this section, we give a comment on Proposition~\ref{prop:genecon}.   
The author could not find any nontrivial example of skew Hadamard difference sets fitting the general construction given in Proposition~\ref{prop:genecon}  other than those obtained from Theorem~\ref{thm:mainskew}. 
Let us consider pure Gauss sums not satisfying  property $(\star)$, e.g., Gauss sums in 
the $r=2$ case such that $(2p_i,f,p)\not\in {\mathcal P}$ for each $i=1,2$ and  $(2p_1p_2,f,p)\in {\mathcal P}$. Note that  the pure Gauss sums in Proposition~\ref{prop:88} (2) belong to this class. 
Then, the conditions (1) and (2) in Proposition~\ref{prop:genecon}  are equivalent to that 
$\sum_{i\in I}\zeta_{t}^i=0$ for any $t\in \{p_1,p_2,p_1p_2,2p_1,2p_2\}$. 
By identifying the subset $I$ with the polynomial $f(x)=\sum_{i\in I}x^i\,(\mod{x^N-1})$,  the condition above is equivalent to   
\begin{equation}\label{eq:divis}
f(x)\equiv 0\,(\mod{\Phi_t}), \quad\, \, \, \forall t\in\{p_1,p_2,p_1p_2,2p_1,2p_2\}, 
\end{equation}
where $\Phi_{t}$ is the $t$th cyclotomic polynomial. The problem is whether there is a polynomial $f(x)\,(\mod{x^N-1})$ with coefficients from 
$\{0,1\}$ and with exactly $p_1p_2$ nonzero coefficients such that 
$f(x)\not\equiv 0\,(\mod{\Phi_{2p_1p_2}})$ and \eqref{eq:divis} is satisfied. 
For example, we checked by a computer that there is no such $f(x)$ for  $(p_1,p_2)=(3,5),(3,7)$. This problem remains open in general, which is difficult but interesting besides evaluating Gauss sums. 

\newpage
\section*{Appendix} 
In this appendix, we list all $(N,\overline{p})\in {\mathcal P}^\ast_f$ for $N\le 5000$ and odd $f$ in Tables~\ref{Tab0} and \ref{Tab01}. 
{\tiny 
\begin{table}[h]
\caption{$(N,\overline{p})\in {\mathcal P}^\ast_f$ for $N\le 5000$ and odd $f$}
\label{Tab0}
\begin{tabular}{|c|c|}
\hline
$[N,f,\overline{p}]$&Ref.\\
\hline 
\hline 
$[ 14, 3, 9 ]$ &  Prop.~\ref{prop:22} \\
$[ 42, 3, 25 ]$ &  Prop.~\ref{prop:44} (2) \\
$[ 46, 11, 3 ]$ &  Prop.~\ref{prop:22} \\
$[ 62, 15, 7 ]$ & Prop.~\ref{prop:22}  \\
$[ 62, 5, 33 ]$ &  Prop.~\ref{prop:66} \\
$[ 78, 3, 55 ]$ & Prop.~\ref{prop:88} (2)  \\ 
$[ 94, 23, 3 ]$ & Prop.~\ref{prop:22}  \\
$[ 98, 21, 9 ]$ & Prop.~\ref{prop:22}  \\
$[ 110, 5, 31 ]$ & Prop.~\ref{prop:88} (2)  \\ 
$[ 142, 35, 3 ]$ & Prop.~\ref{prop:22}  \\
$[ 146, 9, 37 ]$ & Prop.~\ref{prop:88} (1)   \\
$[ 158, 39, 5 ]$ & Prop.~\ref{prop:22}   \\
$[ 178, 11, 39 ]$ & Prop.~\ref{prop:88} (1)   \\
$[ 186, 15, 7 ]$ & Prop.~\ref{prop:44} (2) \\
$[ 206, 51, 7 ]$ & Prop.~\ref{prop:22} \\
$[ 254, 63, 9 ]$ & Prop.~\ref{prop:22}   \\
$[ 254, 21, 25 ]$ & Prop.~\ref{prop:66}   \\
$[ 254, 7, 129 ]$ & ${\mathcal P}^{(2)}$  \\
$[ 294, 21, 25 ]$ & Prop.~\ref{prop:44} (2)  \\
$[ 302, 75, 5 ]$ & Prop.~\ref{prop:22}    \\
$[ 302, 15, 85 ]$ & ${\mathcal P}^{(2)}$  \\
$[ 322, 33, 9 ]$ & Prop.~\ref{prop:44} i)  \\
$[ 334, 83, 3 ]$ & Prop.~\ref{prop:22}    \\
$[ 382, 95, 3 ]$ & Prop.~\ref{prop:22}   \\
$[ 398, 99, 7 ]$ & Prop.~\ref{prop:22}   \\
$[ 434, 15, 39 ]$ & ${\mathcal P}^{(2)}$  \\
$[ 446, 37, 7 ]$ & Prop.~\ref{prop:66}   \\
$[ 446, 111, 9 ]$ & Prop.~\ref{prop:22}   \\
$[ 462, 15, 25 ]$ & Prop.~\ref{prop:88} (3)-iii)  \\ 
$[ 466, 29, 19 ]$ & Prop.~\ref{prop:88} (1) \\
$[ 474, 39, 13 ]$ & Prop.~\ref{prop:44} (2)  \\
$[ 478, 119, 3 ]$ & Prop.~\ref{prop:22}  \\
$[ 506, 55, 3 ]$ & Prop.~\ref{prop:44} (2)  \\
$[ 526, 131, 3 ]$ & Prop.~\ref{prop:22}  \\
$[ 542, 135, 7 ]$ & Prop.~\ref{prop:22}   \\
$[ 618, 51, 7 ]$ &  Prop.~\ref{prop:44} (2)  \\
$[ 622, 155, 3 ]$ &  Prop.~\ref{prop:22}  \\
$[ 654, 27, 7 ]$ & Prop.~\ref{prop:88} (2)  \\
$[ 658, 69, 9 ]$ & Prop.~\ref{prop:44} (1)  \\
$[ 674, 21, 13 ]$ & ${\mathcal P}^{(2)}$ \\ 
$[ 686, 147, 9 ]$ & Prop.~\ref{prop:22}  \\
$[ 718, 179, 3 ]$ & Prop.~\ref{prop:22}  \\
$[ 734, 183, 13 ]$ & Prop.~\ref{prop:22}  \\
$[ 762, 63, 13 ]$ & Prop.~\ref{prop:44} (2)  \\
$[ 766, 191, 3 ]$ &  Prop.~\ref{prop:22} \\
\hline
\end{tabular}
\begin{tabular}{|c|c|}
\hline
$[N,f,\overline{p}]$&Ref.\\
\hline
\hline 
$[ 826, 87, 9 ]$ & Prop.~\ref{prop:44} (2)  \\
$[ 862, 43, 3 ]$ & ${\mathcal P}^{(2)}$  \\
$[ 862, 215, 5 ]$ & Prop.~\ref{prop:22}  \\
$[ 874, 99, 9 ]$ & Prop.~\ref{prop:44} (2)   \\
$[ 878, 219, 5 ]$ &  Prop.~\ref{prop:22}  \\    
$[ 878, 73, 7 ]$ &  Prop.~\ref{prop:66}  \\
$[ 906, 75, 25 ]$ & Prop.~\ref{prop:44} (2)  \\
$[ 926, 231, 9 ]$ & Prop.~\ref{prop:22}  \\
$[ 958, 239, 3 ]$ & Prop.~\ref{prop:22}   \\
$[ 974, 243, 9 ]$ &  Prop.~\ref{prop:22} \\
$[ 994, 105, 9 ]$ & Prop.~\ref{prop:44} (1)  \\
$[ 1006, 251, 3 ]$ & Prop.~\ref{prop:22}   \\
$[ 1014, 39, 55 ]$ & Prop.~\ref{prop:88} (2)  \\ 
$[ 1022, 9, 513 ]$ & ${\mathcal P}^{(2)}$   \\
$[ 1034, 115, 3 ]$ & Prop.~\ref{prop:44} (2)  \\ 
$[ 1058, 253, 3 ]$ & Prop.~\ref{prop:22}  \\
$[ 1086, 45, 13 ]$ & Prop.~\ref{prop:88} (2)  \\ 
$[ 1106, 39, 11 ]$ & ${\mathcal P}^{(2)}$   \\
$[ 1162, 123, 9 ]$ & Prop.~\ref{prop:44} (2)  \\
$[ 1194, 99, 7 ]$ & Prop.~\ref{prop:44} (2)  \\
$[ 1198, 299, 3 ]$ & Prop.~\ref{prop:22}  \\
$[ 1202, 75, 3 ]$ & Prop.~\ref{prop:88} (1)  \\
$[ 1202, 25, 27 ]$ & ${\mathcal P}^{(2)}$  \\
$[ 1210, 55, 31 ]$ & Prop.~\ref{prop:88} (2)  \\ 
$[ 1214, 303, 9 ]$ &  Prop.~\ref{prop:22} \\
$[ 1246, 33, 39 ]$ &  ${\mathcal P}^{(2)}$ \\
$[ 1262, 315, 9 ]$ & Prop.~\ref{prop:22}   \\
$[ 1262, 45, 47 ]$ & ${\mathcal P}^{(2)}$  \\
$[ 1294, 323, 3 ]$ & Prop.~\ref{prop:22} \\
$[ 1310, 65, 11 ]$ &  Prop.~\ref{prop:88} (2)   \\ 
$[ 1338, 111, 19 ]$ &  Prop.~\ref{prop:44} (2)  \\
$[ 1374, 57, 19 ]$ &  Prop.~\ref{prop:88} (2)   \\  
$[ 1426, 165, 9 ]$ &  Prop.~\ref{prop:44} (1)   \\
$[ 1426, 55, 35 ]$ &  ${\mathcal P}^{(2)}$  \\
$[ 1438, 359, 3 ]$ &  Prop.~\ref{prop:22} \\  
$[ 1442, 51, 121 ]$ &  ${\mathcal P}^{(2)}$  \\
$[ 1454, 363, 7 ]$ &  Prop.~\ref{prop:22}   \\
$[ 1454, 121, 9 ]$ &  Prop.~\ref{prop:66}   \\
$[ 1486, 371, 3 ]$ &  Prop.~\ref{prop:22}  \\
$[ 1502, 375, 5 ]$ &  Prop.~\ref{prop:22}  \\
$[ 1626, 135, 7 ]$ &  Prop.~\ref{prop:44} (2)  \\
$[ 1646, 411, 9 ]$ &  Prop.~\ref{prop:22}  \\
$[ 1662, 69, 49 ]$ &  Prop.~\ref{prop:88} (2)   \\ 
$[ 1678, 419, 3 ]$ &  Prop.~\ref{prop:22}     \\
$[ 1726, 431, 3 ]$ &  Prop.~\ref{prop:22}   \\
\hline
\end{tabular}
\begin{tabular}{|c|c|}
\hline
$[N,f,\overline{p}]$&Ref.\\
\hline
\hline 
$[ 1762, 55, 21 ]$ & ${\mathcal P}^{(2)}$  \\
$[ 1774, 443, 3 ]$ &   Prop.~\ref{prop:22}   \\
$[ 1778, 21, 135 ]$ &  ${\mathcal P}^{(2)}$ \\
$[ 1786, 207, 9 ]$ &  Prop.~\ref{prop:44} (2)  \\
$[ 1822, 455, 3 ]$ &  Prop.~\ref{prop:22}   \\
$[ 1822, 91, 15 ]$ & ${\mathcal P}^{(2)}$  \\
$[ 1834, 195, 9 ]$ &  Prop.~\ref{prop:44} (2) \\
$[ 1838, 459, 5 ]$ &  Prop.~\ref{prop:22}  \\
$[ 1838, 153, 9 ]$ &  Prop.~\ref{prop:66}  \\
$[ 1874, 117, 9 ]$ & Prop.~\ref{prop:88} (1)    \\
$[ 1922, 465, 7 ]$ &  Prop.~\ref{prop:22}    \\
$[ 1922, 155, 33 ]$ &  Prop.~\ref{prop:66}  \\
$[ 1934, 483, 21 ]$ &  Prop.~\ref{prop:22}  \\
$[ 1966, 491, 3 ]$ &   Prop.~\ref{prop:22} \\
$[ 1978, 231, 9 ]$ &   Prop.~\ref{prop:44} (2) \\
$[ 1982, 495, 5 ]$ &  Prop.~\ref{prop:22}   \\
$[ 2058, 147, 25 ]$ &   Prop.~\ref{prop:44} (2)  \\
$[ 2062, 515, 3 ]$ & Prop.~\ref{prop:22}  \\
$[ 2078, 519, 7 ]$ &  Prop.~\ref{prop:22}  \\
$[ 2110, 105, 51 ]$ &  Prop.~\ref{prop:88} (2)   \\ 
$[ 2114, 75, 25 ]$ &  ${\mathcal P}^{(2)}$  \\
$[ 2114, 15, 529 ]$ &  ${\mathcal P}^{(2)}$  \\
$[ 2126, 531, 9 ]$ &  Prop.~\ref{prop:22}   \\
$[ 2162, 253, 3 ]$ &  Prop.~\ref{prop:44} (1)   \\
$[ 2174, 543, 9 ]$ &  Prop.~\ref{prop:22}    \\
$[ 2202, 183, 13 ]$ &  Prop.~\ref{prop:44} (2)   \\
$[ 2206, 551, 3 ]$ & Prop.~\ref{prop:22}  \\
$[ 2206, 29, 69 ]$ & ${\mathcal P}^{(2)}$  \\
$[ 2254, 231, 9 ]$ & Prop.~\ref{prop:44} (1)    \\
$[ 2266, 255, 15 ]$ &  Prop.~\ref{prop:44} (2)   \\
$[ 2302, 575, 3 ]$ &  Prop.~\ref{prop:22}   \\
$[ 2338, 249, 9 ]$ &  Prop.~\ref{prop:44} (1) \\
$[ 2446, 611, 7 ]$ &  Prop.~\ref{prop:22} \\
$[ 2462, 615, 5 ]$ &  Prop.~\ref{prop:22} \\
$[ 2478, 87, 25 ]$ &  Prop.~\ref{prop:88} (3)-ii) \\  
$[ 2510, 125, 21 ]$ &   Prop.~\ref{prop:88} (2) \\ 
$[ 2526, 105, 25 ]$ & Prop.~\ref{prop:88} (2)  \\ 
$[ 2558, 639, 5 ]$ &  Prop.~\ref{prop:22}  \\
$[ 2578, 161, 29 ]$ &  Prop.~\ref{prop:88} (1) \\
$[ 2606, 651, 23 ]$ & Prop.~\ref{prop:22} \\
$[ 2622, 99, 25 ]$ &  Prop.~\ref{prop:88} (3)-iii)  \\  
$[ 2634, 219, 13 ]$ &  Prop.~\ref{prop:44} (2)\\
$[ 2638, 659, 3 ]$ & Prop.~\ref{prop:22} \\
$[ 2654, 663, 9 ]$ &  Prop.~\ref{prop:22}\\
$[ 2654, 221, 43 ]$ &  Prop.~\ref{prop:66} \\
\hline
\end{tabular}
\end{table}}

\newpage 

{\tiny
\begin{table}[h]
\caption{$(N,\overline{p})\in {\mathcal P}^\ast_f$ for $N\le 5000$ and odd $f$}
\label{Tab01}
\begin{tabular}{|c|c|}
\hline
$[N,f,\overline{p}]$&Ref.\\
\hline
\hline 
$[ 2674, 285, 9 ]$ &  Prop.~\ref{prop:44} (2)\\
$[ 2734, 683, 3 ]$ &  Prop.~\ref{prop:22} \\
$[ 2778, 231, 25 ]$ &  Prop.~\ref{prop:44} (2) \\
$[ 2782, 159, 3 ]$ &  Prop.~\ref{prop:88} (2) \\ 
$[ 2786, 99, 23 ]$ &   ${\mathcal P}^{(2)}$ \\
$[ 2798, 699, 5 ]$ & Prop.~\ref{prop:22}  \\
$[ 2798, 233, 9 ]$ &  Prop.~\ref{prop:66}  \\
$[ 2846, 711, 9 ]$ &  Prop.~\ref{prop:22} \\
$[ 2846, 237, 23 ]$ &  Prop.~\ref{prop:66} \\
$[ 2866, 179, 15 ]$ &  Prop.~\ref{prop:88} (1) \\
$[ 2878, 719, 3 ]$ &  Prop.~\ref{prop:22} \\
$[ 2894, 723, 9 ]$ &  Prop.~\ref{prop:22} \\
$[ 2914, 345, 7 ]$ &  Prop.~\ref{prop:44} (1) \\
$[ 2914, 115, 97 ]$ &   ${\mathcal P}^{(2)}$ \\
$[ 2922, 243, 31 ]$ &  Prop.~\ref{prop:44} (2) \\
$[ 2942, 735, 5 ]$ &  Prop.~\ref{prop:22} \\
$[ 2942, 245, 19 ]$ &  Prop.~\ref{prop:66} \\
$[ 2974, 743, 3 ]$ &  Prop.~\ref{prop:22} \\
$[ 3022, 755, 5 ]$ & Prop.~\ref{prop:22} \\
$[ 3038, 105, 39 ]$ &   Prop.~\ref{prop:44} (1) \\
$[ 3086, 771, 13 ]$ &  Prop.~\ref{prop:22}  \\
$[ 3118, 779, 3 ]$ &  Prop.~\ref{prop:22}  \\
$[ 3122, 111, 289 ]$ &   ${\mathcal P}^{(2)}$ \\
$[ 3134, 783, 7 ]$ &  Prop.~\ref{prop:22}  \\
$[ 3166, 791, 11 ]$ &  Prop.~\ref{prop:22}  \\
$[ 3178, 339, 9 ]$ &  Prop.~\ref{prop:44} (2) \\
$[ 3214, 803, 3 ]$ &  Prop.~\ref{prop:22}  \\
$[ 3218, 201, 11 ]$ & Prop.~\ref{prop:88} (1) \\
$[ 3234, 105, 25 ]$ &  Prop.~\ref{prop:88} (3)-iii) \\  
$[ 3246, 135, 25 ]$ & Prop.~\ref{prop:88} (2)   \\ 
$[ 3262, 87, 23 ]$ &  ${\mathcal P}^{(2)}$  \\
$[ 3266, 385, 3 ]$ &Prop.~\ref{prop:44} (1)  \\
$[ 3310, 165, 21 ]$ & Prop.~\ref{prop:88} (2)   \\ 
$[ 3326, 831, 9 ]$ & Prop.~\ref{prop:22} \\
$[ 3346, 357, 9 ]$ & Prop.~\ref{prop:44} (1) \\
$[ 3358, 99, 55 ]$ &  ${\mathcal P}^{(2)}$\\
$[ 3406, 195, 3 ]$ & Prop.~\ref{prop:88} (2) \\ 
$[ 3422, 203, 7 ]$ & Prop.~\ref{prop:88} (2) \\ 
$[ 3442, 215, 17 ]$ & Prop.~\ref{prop:88} (1) \\
\hline
\end{tabular}
\begin{tabular}{|c|c|}
\hline
$[N,f,\overline{p}]$&Ref.\\
\hline
\hline 
$[ 3454, 195, 9 ]$ & Prop.~\ref{prop:88} (2)  \\ 
$[ 3486, 123, 25 ]$ & Prop.~\ref{prop:88} (3)-ii) \\ 
$[ 3514, 375, 9 ]$ &  Prop.~\ref{prop:44} (2) \\
$[ 3518, 879, 11 ]$ & Prop.~\ref{prop:22}\\
$[ 3566, 891, 7 ]$ & Prop.~\ref{prop:22} \\
$[ 3602, 225, 9 ]$ & Prop.~\ref{prop:88} (1)\\
$[ 3602, 75, 21 ]$ &  ${\mathcal P}^{(2)}$ \\
$[ 3602, 25, 175 ]$ &  ${\mathcal P}^{(2)}$\\
$[ 3634, 429, 9 ]$ & Prop.~\ref{prop:44} (1) \\
$[ 3642, 303, 13 ]$ & Prop.~\ref{prop:44} (2) \\
$[ 3646, 911, 3 ]$ & Prop.~\ref{prop:22}  \\
$[ 3662, 305, 5 ]$ & Prop.~\ref{prop:66} \\
$[ 3662, 915, 9 ]$ & Prop.~\ref{prop:22} \\
$[ 3682, 393, 9 ]$ & Prop.~\ref{prop:44} (1) \\
$[ 3694, 923, 3 ]$ & Prop.~\ref{prop:22} \\
$[ 3742, 935, 5 ]$ & Prop.~\ref{prop:22}  \\
$[ 3758, 939, 5 ]$ &  Prop.~\ref{prop:22} \\
$[ 3786, 315, 31 ]$ & Prop.~\ref{prop:44} (2) \\
$[ 3794, 135, 37 ]$ & ${\mathcal P}^{(2)}$ \\
$[ 3818, 451, 3 ]$ & Prop.~\ref{prop:44} (2) \\
$[ 3826, 239, 17 ]$ & Prop.~\ref{prop:88} (1) \\
$[ 3838, 225, 5 ]$ & Prop.~\ref{prop:88} (2) \\ 
$[ 3902, 975, 5 ]$ & Prop.~\ref{prop:22} \\
$[ 3998, 999, 5 ]$ & Prop.~\ref{prop:22} \\
$[ 3998, 333, 13 ]$ & Prop.~\ref{prop:66} \\
$[ 4042, 483, 9 ]$ & Prop.~\ref{prop:44} (2)  \\
$[ 4042, 161, 21 ]$ & Exception \\
$[ 4078, 1019, 3 ]$ & Prop.~\ref{prop:22}  \\
$[ 4094, 11, 2049 ]$ & ${\mathcal P}^{(2)}$ \\
$[ 4126, 1031, 3 ]$ & Prop.~\ref{prop:22}  \\
$[ 4174, 1043, 3 ]$ & Prop.~\ref{prop:22} \\
$[ 4178, 261, 15 ]$ & Prop.~\ref{prop:88} (1) \\
$[ 4178, 87, 85 ]$ & ${\mathcal P}^{(2)}$ \\
$[ 4178, 29, 457 ]$ & ${\mathcal P}^{(2)}$ \\
$[ 4222, 1055, 5 ]$ & Prop.~\ref{prop:22} \\
$[ 4254, 177, 7 ]$ &  Prop.~\ref{prop:88} (2)  \\ 
$[ 4286, 1071, 9 ]$ &  Prop.~\ref{prop:22} \\
$[ 4286, 357, 15 ]$ &  Prop.~\ref{prop:66} \\
$[ 4286, 51, 67 ]$ & ${\mathcal P}^{(2)}$   \\
\hline
\end{tabular}
\begin{tabular}{|c|c|}
\hline
$[N,f,\overline{p}]$&Ref.\\
\hline
\hline 
$[ 4286, 153, 121 ]$ & ${\mathcal P}^{(2)}$ \\
$[ 4354, 465, 9 ]$ & Prop.~\ref{prop:44} (1) \\
$[ 4362, 363, 7 ]$ & Prop.~\ref{prop:44} (2)  \\
$[ 4378, 495, 9 ]$ & Prop.~\ref{prop:44} (2) \\
$[ 4398, 183, 31 ]$ & Prop.~\ref{prop:88} (2) \\ 
$[ 4402, 105, 19 ]$ & ${\mathcal P}^{(2)}$ \\
$[ 4402, 35, 225 ]$ &  ${\mathcal P}^{(2)}$ \\
$[ 4414, 1103, 3 ]$ & Prop.~\ref{prop:22}  \\
$[ 4418, 1081, 3 ]$ & Prop.~\ref{prop:22}  \\
$[ 4478, 1119, 9 ]$ & Prop.~\ref{prop:22} \\
$[ 4494, 159, 25 ]$ & Prop.~\ref{prop:88} (3)-iii) \\ 
$[ 4506, 375, 13 ]$ &  Prop.~\ref{prop:44} (2) \\
$[ 4526, 45, 221 ]$ & ${\mathcal P}^{(2)}$  \\
$[ 4542, 189, 7 ]$ & Prop.~\ref{prop:88} (2) \\ 
$[ 4574, 1143, 7 ]$ & Prop.~\ref{prop:22}  \\
$[ 4574, 381, 9 ]$ & Prop.~\ref{prop:66}  \\
$[ 4606, 483, 9 ]$ & Prop.~\ref{prop:44} (1) \\
$[ 4622, 1155, 9 ]$ & Prop.~\ref{prop:22} \\
$[ 4702, 1175, 3 ]$ & Prop.~\ref{prop:22} \\
$[ 4702, 235, 11 ]$ & ${\mathcal P}^{(2)}$ \\
$[ 4702, 47, 15 ]$ & ${\mathcal P}^{(2)}$ \\
$[ 4718, 21, 1115 ]$ &${\mathcal P}^{(2)}$ \\
$[ 4738, 561, 25 ]$ & Prop.~\ref{prop:44} (1)  \\
$[ 4766, 397, 7 ]$ & Prop.~\ref{prop:66} \\
$[ 4766, 1191, 13 ]$ & Prop.~\ref{prop:22} \\
$[ 4798, 1199, 3 ]$ & Prop.~\ref{prop:22} \\
$[ 4802, 1029, 9 ]$ & Prop.~\ref{prop:22} \\
$[ 4814, 287, 7 ]$ & Prop.~\ref{prop:88} (2)  \\ 
$[ 4846, 1211, 3 ]$ & Prop.~\ref{prop:22} \\
$[ 4882, 305, 5 ]$ & Prop.~\ref{prop:88} (1) \\
$[ 4894, 1223, 3 ]$ & Prop.~\ref{prop:22} \\
$[ 4898, 195, 95 ]$ &  ${\mathcal P}^{(2)}$ \\
$[ 4906, 555, 9 ]$ & Prop.~\ref{prop:44} (2)  \\
$[ 4910, 245, 11 ]$ & Prop.~\ref{prop:88} (2) \\ 
$[ 4922, 583, 3 ]$ & Prop.~\ref{prop:44} (2) \\
$[ 4938, 411, 13 ]$ & Prop.~\ref{prop:44} (2) \\
$[ 4974, 207, 55 ]$ & Prop.~\ref{prop:88} (2)  \\
 &  \\
 &   \\
\hline
\end{tabular}
\end{table}}

\end{document}